\documentclass[11pt]{article}
\usepackage[T1]{fontenc}
\usepackage{latexsym,amssymb,amsmath,amsfonts,amsthm}
\usepackage{graphicx}
\usepackage{epstopdf}
\usepackage{float}
\usepackage{subfigure}
\usepackage{caption}
\usepackage{amsmath}
\numberwithin{figure}{section}
\numberwithin{table}{section}
\usepackage{placeins}
\usepackage{color}
\usepackage{indentfirst}
\numberwithin{equation}{section}
\newcommand{\R}{{\mathbb R}}

\newcommand{\N}{{\mathbb N}}
\newcommand{\C}{{\mathbb C}}
\newcommand{\s}{{\mathbb S}}

\newcommand{\be}{\begin{eqnarray}}
\newcommand{\ben}{\begin{eqnarray*}}
\newcommand{\en}{\end{eqnarray}}
\newcommand{\enn}{\end{eqnarray*}}

\newtheorem{theorem}{Theorem}[section]
\newtheorem{lemma}[theorem]{Lemma}
\newtheorem{corollary}[theorem]{Corollary}
\newtheorem{definition}[theorem]{Definition}
\newtheorem{remark}[theorem]{Remark}

\definecolor{rot}{rgb}{1.000,0.000,0.000}
\definecolor{blau}{rgb}{0,0,1}

\begin{document}
	\captionsetup[figure]{labelfont={bf},name={Fig.},labelsep=period}
	
	\title{Factorization method with one plane wave: from model-driven and data-driven perspectives}

	\author{Guanqiu Ma\footnotemark[1]\  \and Guanghui Hu\footnotemark[2]
}
	
	\date{}
	\maketitle
	
	\renewcommand{\thefootnote}{\fnsymbol{footnote}}
	\footnotetext[1]{Department of Applied Mathematics, Beijing Computational Science Research Center, Beijing 100193, P. R. China (guanqium@csrc.ac.cn).}
\footnotetext[2]{School of Mathematical Sciences and LPMC, Nankai Universtiy, Tianjin 300071, China (ghhu@nankai.edu.cn, corresponding author).}

	\renewcommand{\thefootnote}{\arabic{footnote}}

	\begin{abstract}
The factorization method by Kirsch (1998) provides a necessary and sufficient condition for characterizing the shape and position of an unknown scatterer by using far-field patterns of infinitely many time-harmonic plane waves at a fixed frequency. This paper is concerned with the factorization method with a single far-field pattern to recover a convex polygonal scatterer/source. Its one-wave version relies on the absence of analytical continuation of the scattered/radiated wave-fields in corner domains.
 It can be regarded as a domain-defined sampling method and does not require forward solvers. In this paper we provide a rigorous mathematical justification of the one-wave factorization method and present some preliminary numerical examples. In particular, the proposed scheme 
 can be interpreted as a model-driven and data-driven method, because it essentially
 depends on the scattering model and a priori given \emph{sample data}.

\vspace{.2in} {\bf Keywords}: factorization method,	inverse scattering, inverse source problem, single far-field pattern, polygonal scatterers, corner scattering.
	\end{abstract}
	\section{Introduction}
\begin{figure}[h]
 	\centering
 	\subfigure
 	{
 		\begin{minipage}[t]{0.4\linewidth}
 			\centering
 			\includegraphics[scale=0.3]{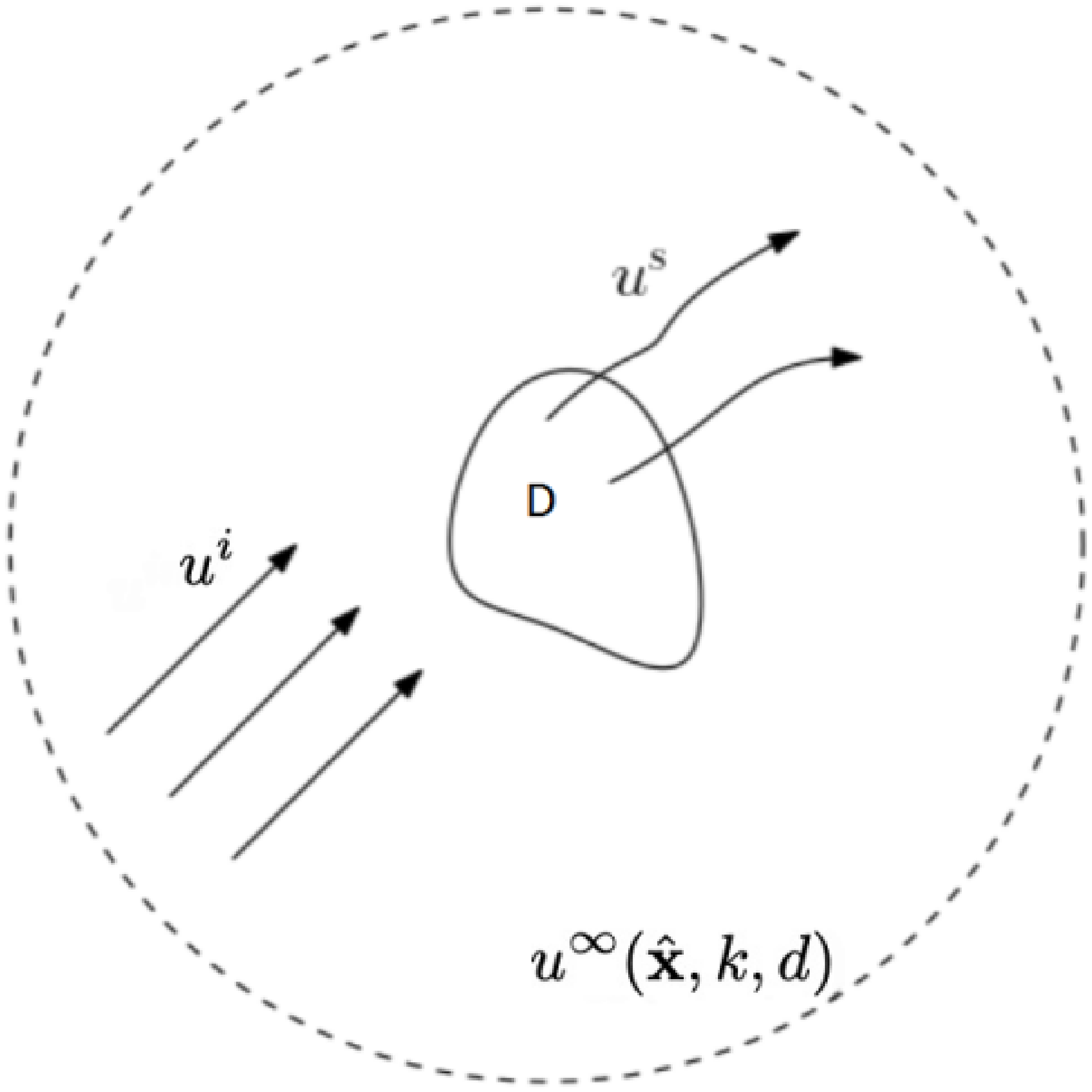}
 		\end{minipage}
 	}
 	\subfigure
 	{
 		\begin{minipage}[t]{0.4\linewidth}
 			\centering
 			\includegraphics[scale=0.3]{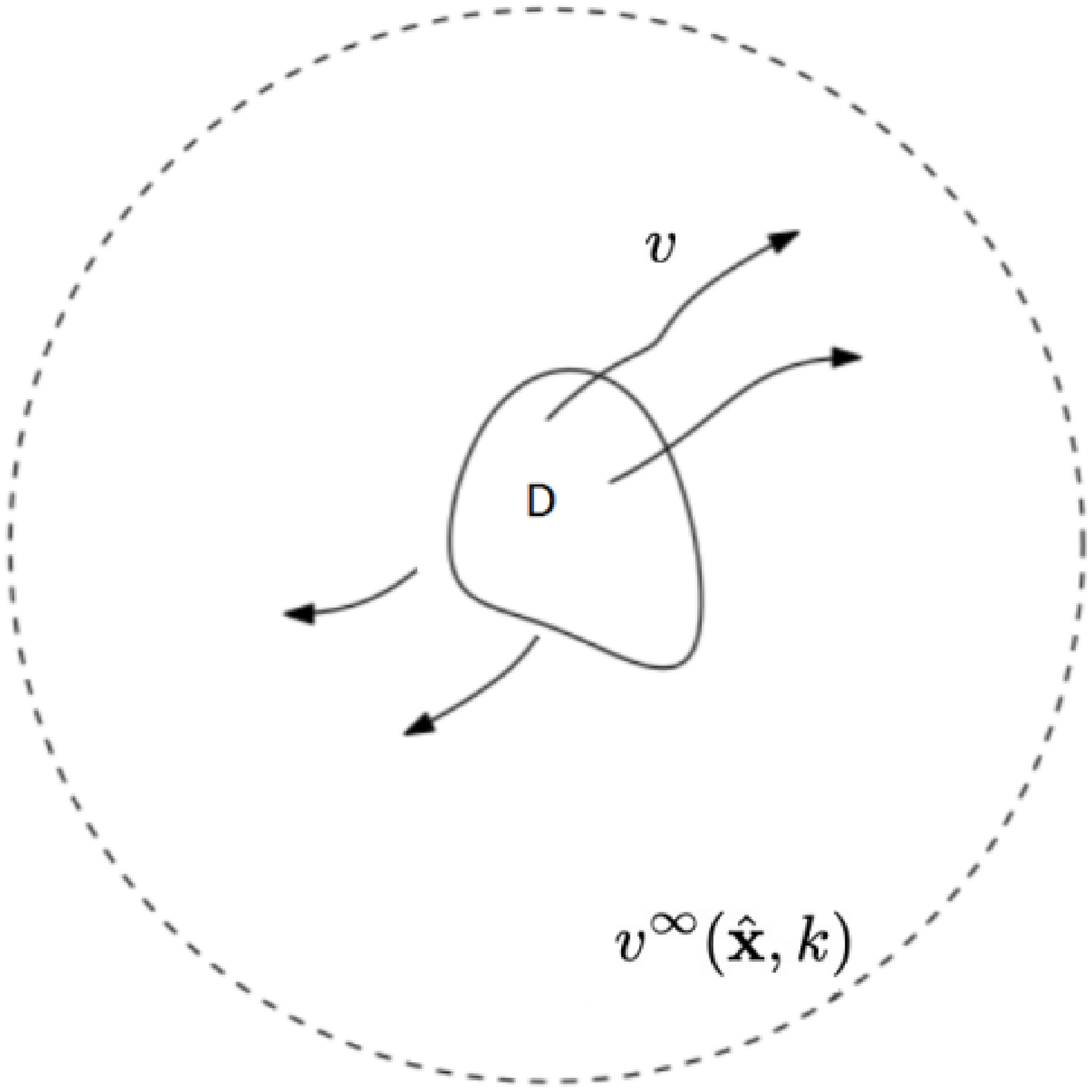}
 		\end{minipage}
 	}
 \caption{Illustration of an obstacle scattering problem (left) and a source radiating problem (right).}
 \label{scattering}
 \end{figure}

The primary goal of inverse scattering theory is to extract information about unknown objects from the wave-fields measured far way from the target. Time-harmonic inverse scattering is widely used in deep-sea exploration, geological exploration, medical imaging, non-destructive test and other fields. To describe the two-dimensional model for inverse obstacle scattering problems, 
we consider the propagation of a time-harmonic incident field $u^{i}$ in a homogeneous and isotropic background medium governed by the Helmholtz equation
	\begin{equation}
		\Delta u^{i} + k^{2}u^{i} = 0 \hspace{.2 cm} \text{in} \hspace{.2 cm} \mathbb{R}^{2},
	\end{equation}
	where $k>0$ is the wavenumber. Assume that a plane wave  $u^{i}=e^{ikx\cdot d}$ with the direction $d=(\cos\theta,\sin\theta)$ ($\theta\in[0,2\pi]$) is incident onto a sound-soft scatterer $D\subset \mathbb{R}^2$; see Figure \ref{scattering} (left). The scatterer $D$ is supposed to occupy a bounded Lipschitz domain such that its exterior $\mathbb{R}^{2} \backslash \overline{D}$ is connected. 
	The scattered field $u^{s}$ is also governed by the Helmholtz equation
	\begin{equation}
		\Delta u^{s} + k^{2}u^{s} = 0 \hspace{.2 cm} \text{in} \hspace{.2 cm} \mathbb{R}^{2}\backslash \overline{D},
		\label{Helmholtz_obstacle_eq}
	\end{equation}
	and satisfies the Dirichlet boundary condition
	\begin{equation}
		u^{s} =  -u^{i} \hspace{.2 cm} \text{on} \hspace{.2 cm} \partial D
		\label{DiriCond_Helmholtz_eq}
	\end{equation}
	together with the outgoing Sommerfeld radiation condition
	\begin{equation}
		\lim_{r\rightarrow \infty} \sqrt{r} \left( \frac{\partial u^{s}}{\partial r} - iku^{s} \right) = 0, \hspace{.2 cm} r = |x|,
		\label{SommerfeldCond}
	\end{equation}
uniformly in all directions $\hat{x} = x/|x|\in \mathbb{S}:=\{x: |x|=1\}$, $x \in \mathbb{R}^{2} \backslash \overline{D}$. The Sommerfeld radiation condition of $u^{s}$  leads to an asymptotic behavior of $u^s$ in the form
    \begin{equation}
    u^{s}(x) = \frac{ e^{\mathrm{i}k|x|} }{ \sqrt{|x|} } \left\{   u^{\infty}(\hat{x}) + \mathcal{O}\left(\frac{ 1 }{ \sqrt{|x|} }\right)  \right\}, \hspace{.2 cm} |x| \rightarrow \infty,
    \label{asymptotic_behaviour}
    \end{equation}
    where $u^{\infty}(\hat{x})=u^{\infty}(\hat{x},k,d)$ is called the far-field pattern at the observation direction $\hat{x}\in\mathbb{S}$. The total field $u$ is defined as $u:= u^{i} + u^{s}$ in $\mathbb{R}^{2}\backslash \overline{D}$. Using variational or integral equation method, it is well known that the system \eqref{Helmholtz_obstacle_eq}-\eqref{asymptotic_behaviour} always admits a unique solution $u^s\in H^1_{loc}(\R^2\backslash\overline{D})$; see e.g. \cite{Cakoni,kress.1998,Kirsch.2008,NP2013,Potthast}.
    The goal of inverse obstacle scattering is to recover $D$ from far-field patterns incited by one or many incoming waves.

    For the model of inverse source problems, we consider the radiated field $v$  governed by the inhomogeneous Helmholtz equation (see Figure \ref{scattering} (right))
    	\begin{equation}
    	\Delta v + k^{2}v = f \hspace{.2 cm} \text{in} \hspace{.2 cm} \mathbb{R}^{2},
    	\label{Helmholtz_source_eq}
    \end{equation}
    where the source term $f\in L^2_{loc}(\R^2)$ is supposed to be compactly supported on $\overline{D}$ (that is, $D=\mbox{supp}(f)$). Here $v$ is required to fulfill the Sommerfeld radiation condition \eqref{SommerfeldCond}. The inverse source problem is to recover $f$ or its support $D$ from the far-field pattern $v^\infty$ of $v$.

In this paper we are interested in non-iterative approaches for recovering $\partial D$ from a single far-field pattern. Such inverse problems are well-known to be nonlinear and ill-posed. In comparison with the optimization-based iterative schemes, sampling methods (which are also called qualitative methods in literature \cite{Cakoni}) have attracted much attention over the last twenty years, since they do not need any forward solver and initial approximation of the target. Basically there exist two kinds of sampling methods: multi-wave and one-wave sampling methods.  The multi-wave sampling methods do not need a priori information on physical and geometrical properties of the scatterer, but usually require the knowledge of far-field patterns for a large number of incident waves.
They consists of both
point-wisely defined and domain-defined inversion schemes.
 The first ones are usually based on designing an appropriate point-wisely defined indicator function which decides on whether a sampling point lies inside or outside of the target. Here we give an incomplete list of such methods: linear sampling method \cite{CK1996,Cakoni}, factorization method \cite{Kirsch98,Kirsch.2008}, singular source method \cite{Potthast},
 orthogonal/direct sampling method \cite{P2010, G2008, IJZ} and the frequency-domain reverse time migration method \cite{CCH}.
 In particular, the factorization method by Kirsch (1998), which has been used in a variety of inverse problems, provides a necessary and sufficient condition for characterizing the shape and position of an unknown scatterer by using multi-static far-field patterns at a fixed frequency. The generalized linear sampling method \cite{glsm} also provides an exact characterization of a scatterer.
 The domain sampling methods are based on choosing an indicator functional which decides on whether a test domain (or a curve) lies inside or outside of the target. Examples of the multi-wave domain sampling methods include, for example,
 range and no-response test \cite{PS2005} and Ikehata's probe method \cite{Ik1998}.

 The one-wave sampling methods are usually designed to test the analytic extensibility of the scattered field; see the monograph \cite[Chapter 15]{NP2013} for detailed discussions. They require only a single far-field pattern or one-pair Cauchy data, but one must pre-assume the absence of an analytical continuation across a general target interface. They are mostly domain-defined sampling methods, for example, range test \cite{rangetest}, no response test \cite{P2003}, enclosure method \cite{Ik1999, I99}  and extended linear sampling method \cite{esm2018}.
 The one-wave range test and no-response test are proven to be dual for both inverse scattering and inverse boundary value problems \cite{NP2013,LNPW}.

 The aim of this paper is to address a framework of the one-wave factorization method, which was earlier proposed in \cite{Elschner_Hu_Elastic} for inverse elastic scattering from rigid polygonal rigid bodies and later discussed in \cite{HL} for inverse acoustic source problems in an inhomogeneous medium. Our arguments are motivated by the existing one-wave sampling methods mentioned above and the recently developed corner scattering theory for justifying the absence of non-scattering energies and non-radiating sources (see \cite{Bl2018,BLS,ElHu2015,hu2018,HSV,HM,KS,LHY,PSV}). The corner scattering theory implies that the wave field cannot be analytically continued across a strongly or weakly singular point lying on the scattering interface.
 The one-wave factorization method leads to a sufficient and necessary condition for imaging any convex penetrable/impenetrable scatterers of polygonal type.
Compared with other one-wave sampling methods, we conclude promising features of the one-wave factorization method as follows. i) The computational criterion involves only inner product calculations and thus looks more straightforward. A new sampling scheme was proposed in this paper by using test disks in two dimensions. Since the number of sampling variables is comparable with the classical linear sampling and factorization methods, the computational cost is not heavier than these multi-wave methods. ii) It is a model-driven and data-driven approach. The one-wave factorization method relies on both the physically scattering model (that is, Helmholtz equation) and the a priori data of some properly chosen test scatterers. In the terminology of learning theory and data sciences, these test scatterers and the associated data are called respectively \emph{samples} and \emph{sample data}.
They are usually given in advance and the sample data can be
calculated off-line before the inversion process. In this paper, we choose sound-soft and impedance disks as test scatterers, because the spectra of the resulting far-field operator take explicit forms. However, there is a variety of choices on the shape and physical properties of test scatterers and also on the type of \emph{sample data}. Finally, it deserves to mentioning that the one-wave factorization method provides a necessary and sufficient condition for recovering convex polygonal scatterers/sources.

    This paper is organized as follows. In Section \ref{sec:2}, we review the multi-wave factorization method for recovering sound-soft and impedance scatterers. In Section \ref{sec:3}, we give a rigorous justification of the one-wave version by combining the classical factorization method and the corner scattering theory. Numerical tests will be performed in Section \ref{sec:4} and concluding remarks will be presented in the final Section \ref{sec:5}.


    \section{Factorization method with infinite plane waves: a model-driven approach}\label{sec:2}
In this section we will briefly review the classical Factorization method \cite{Kirsch98,Kirsch.2008} using the spectra of the far-field operator, which requires measurement data of infinite number of plane waves with distinct directions. It is a typical model-driven scheme, since it depends heavily on the physically scattering model. The resulting computational criterion provides a sufficient and necessary condition for imaging an impenetrable obstacle of sound-soft or impedance type. Below we only review the two-dimensional case. However, the analogue results carry over to three dimensions.
We first recall the Bessel and Neumann functions $J_n$ and $Y_n$, defined as
    \begin{equation}
    J_{n}(t):=\sum_{p=0}^{\infty} \frac{(-1)^{p}}{p !(n+p) !}\left(\frac{t}{2}\right)^{n+2 p}, \\
    Y_{n}(t)=\frac{J_{n}(t) \cos (n \pi)-J_{-n}(t)}{\sin (n \pi)}, \  n \in \mathbb{Z}.
    \end{equation}
The Hankel functions $H_n^{(1)}$ of the first kind of order $n$ are defined as
    \begin{equation}
    H_{n}^{(1)}(t):=J_{n}(t)+i Y_{n}(t).
    \end{equation}
Let $D\subset \R^2$ be the domain occupied by the scatterer.
Recall the single potential operator $S_D$,
    \begin{equation}
    (S_D \psi)(x):=\int_{\partial D} \psi(y) \Phi(x, y) d s(y), \quad x \in \partial D,
    \end{equation}
    where $\Phi(x,y) := \frac{i}{4} H_{0}^{(1)}(k|x-y|), \  x \neq y$ is the fundamental solution of the Helmholtz equation $(\Delta+k^2)u=0$ in $\mathbb{R}^2$. Throughout the paper, the adjoint of an operator will be denoted by $(\cdot)^*$ and the inner product over $L^2(\s)$ by $\langle\cdot,\cdot\rangle_{\s}$.
    \subsection{Factorization method for imaging general scatterers}
    We denote by $u^{\infty}_D (\hat{x}, d)$ the far-field pattern of $u^s$
    to indicate the dependance on the scatterer $D$, which corresponds to the boundary value problems (\ref{Helmholtz_obstacle_eq})-(\ref{SommerfeldCond}) with $u^i=e^{ikx\cdot d}$. Below we state the definition of far-field operator in scattering theory.
	\begin{definition}
	The far-field operator $F_D : L^{2}\left(\mathbb{S}\right) \to L^{2}\left(\mathbb{S}\right)$ corresponding to $D$ is defined by
	$$
	(F_D g)(\hat{x})=\int_\mathbb{S} u^{\infty}_D (\hat{x}, d) g(d) d s(d) \quad \mbox{for all} \ \hat{x} \in \mathbb{S}.
	$$
	\end{definition}
   If $D$ is a sound-soft obstacle, it is well known that $F_D$ is a normal operator. It was proved in \cite[Theorem 1.15]{Kirsch.2008} that the far-field operator $F_D$ can be decomposed into the form
    \begin{equation}\label{fac-D}
          F_D = -G_D S_D^*G_D^*.
    \end{equation}
    Here the data-to-pattern operator $G_D$ : $H^{1/2}(\partial D)\to L^2(\mathbb{S})$ is defined by
    $G_D (f)=v^{\infty}$,
    where $v^{\infty} \in L^2(\mathbb{S})$ is the far-field pattern of the radiation solution $v^s$ to the exterior scattering problem (\ref{Helmholtz_obstacle_eq}) with the boundary data $v^s|_{\partial D}=f \in H^{1/2}(\partial D)$. 
By the Factorization method, the far-field pattern $\phi_z(\hat{x}):=e^{ik\hat{x}\cdot z}$ of the point source wave $x\rightarrow \Phi(x,z)$ belongs to the range of $G_D$ if and only if $z\in D$ (see (see \cite[Theorem 1.12]{Kirsch.2008})). Moreover, the $(F^*F)^{1/4}$-method (see \cite[Theorem 1.24]{Kirsch.2008}) yields the relation $Range(G_D)=Range((F_D^*F_D)^{1/4})$ if $k^2$ is not a Dirichlet eigenvalue of $-\bigtriangleup$ over $D$. Hence, by the Picard theorem, the scatterer $D$ can be characterized by the spectra of $F_D$ as follows.
    \begin{theorem}\label{class_fac}(\cite[Theorem 1.25]{Kirsch.2008})
    	Assume that $k^2$ is not a Dirichlet eigenvalue of $-\triangle$ over $D$. Denote by $(\lambda^{(j)}_{D}, \varphi^{(j)}_{D})$ a spectrum system of the far-field operator $F_D : L^{2}\left(\mathbb{S}\right) \rightarrow L^{2}\left(\mathbb{S}\right)$.
     Then,
       \be\label{FD}
              z \in D  
       \Longleftrightarrow
       I(z) :=\left[ \sum\limits_{j} \frac{\left|\left\langle\phi_{z}, \varphi_{D}^{(j)}\right\rangle_{\mathbb{S}}\right|^{2}}{\left|\lambda_{D}^{(j)}\right|}\right]^{-1}>0,
              \en
	  \end{theorem}
By Theorem \ref{class_fac}, the sign of the indicator function $z\rightarrow I(z)$ can be regarded as the characteristic function of $D$. We note that in (\ref{FD}), $z\in \R^2$ are the sampling variables/points and the spectral data $(\lambda^{(j)}_{D}, \varphi^{(j)}_{D})$ are determined by the far-field patterns $u_D^\infty(\hat{x},d)$ over all observation and incident directions $\hat, d\in \s$.

Now we turn to impenetrable obstacles of impedance type, that is,
\begin{equation}\label{imp-bound}
  	\partial_\nu u+ \eta u = 0 \quad \text{on} \  \partial D,
  \end{equation}
  where $\eta \in L^{\infty}(\partial D)$ is an impedance function satisfying $\operatorname{Im} (\eta) \geq 0$. Denote by $F_{D,imp}: L^{2}\left(\mathbb{S}\right) \to L^{2}\left(\mathbb{S}\right)$ the corresponding far-field operator, and by
$G_{D,imp}: H^{-1/2}(\partial D) \rightarrow L^2(\mathbb{S})$ the
 data-to-pattern operator, that is, 
  \begin{equation}\label{data-pattern-Dimp}
  		G_{D,imp} (f)=v^{\infty},
  \end{equation}
   where $v^{\infty} \in L^2(\mathbb{S})$ is the far-field pattern of the radiation solution $v^s \in H^1_{loc}(\mathbb{R}^2 \backslash \overline{D})$, which solves
\begin{equation}\label{dnp-Dimp-con}
\bigtriangleup v^s + k^2 v^s =0 \quad \text{in} \  \mathbb{R}^2 \backslash \overline{D}, \qquad 	\partial_\nu v^s+ \eta v^s = f \quad \text{on} \  \partial D.
\end{equation}
In the impedance case, the operator $F_{D,imp}$ fails to be normal but can still be factorized
into the form
\begin{equation}\label{fac-Dimp}
F_{D,imp}=-G_{D,imp} T_{D,imp}^{*} G_{D,imp}^{*},
\end{equation}
where $T_{D,imp}: H^{1 / 2}(\partial D) \rightarrow
  H^{-1/2}(\partial D)$ is a Fredholm operator of index zero; see \cite[(2.39)]{Kirsch.2008}.
  Instead of the $(F^*F)^{1/4}$-method for sound-soft obstacles, the $F^\#$-method \cite[Chapter 2.5]{Kirsch.2008} gives an analogous characterization of the impedance obstacle $D$ to the sound-soft case:
   \begin{theorem}\label{class_fac_imp}(\cite[Corollary 2.16]{Kirsch.2008})
  	Assume that $k>0$ is not an eigenvalue of $-\Delta$ over $D$ with respect to the impedance boundary condition with impedance $\lambda$. Then, for any $z \in \mathbb{R}^2$ we have
  	\begin{equation}\label{FI}
  	z \in D  
  	\Longleftrightarrow
  	I_{imp}(z) :=\left[ \sum\limits_{j} \frac{\left|\left\langle\phi_{z}, \varphi^{(j)}_{D}\right\rangle_{\mathbb{S}}\right|^{2}}{\lambda^{(j)}_{D}}\right]^{-1}>0,
  	\end{equation}	
  	where 
  $(\lambda^{(j)}_{D}, \varphi^{(j)}_{D})$ is a spectral system of the positive operator $F_{D,\#} :=
  |{\rm Re} (F_{D,imp})|+|{\rm Im} (F_{D,imp})|: L^{2}\left(\mathbb{S}\right) \rightarrow L^{2}\left(\mathbb{S}\right)$.
  \end{theorem}
  We note that there exists no impedance eigenvalues if ${\rm Im}(\eta)>0$ almost everywhere on an open set of $\partial D$, for instance, $\eta\equiv i\tilde{\eta}$ where $\tilde{\eta}>0$ is a constant.
  In (\ref{FI}), the denominator $\lambda^{(j)}_{D}$ can be equivalently replaced by $|\tilde\lambda^{(j)}_{D}|$ where $\tilde\lambda^{(j)}_{D}$ denote the eigenvalues of $F_{D,imp}$,
  because of the estimate
  \ben
  \frac{1}{\sqrt{2}}\left(|{\rm Re} \tilde\lambda^{(j)}_{D}|+|{\rm Im} \tilde\lambda^{(j)}_{D}|\right)\leq |\tilde\lambda^{(j)}_{D}|\leq |{\rm Re} \tilde\lambda^{(j)}_{D}|+|{\rm Im} \tilde\lambda^{(j)}_{D}|.
  \enn
  \subsection{Explicit examples: Factorization method for imaging disks}\label{FMdisk}
First we give an explicit example for imaging an impedance disk $B_R:=\{x: |x|<R\}$ centered at the origin with radius $R>0$ and the constant impedance coefficient. We suppose that
 the impedance function $\eta=i\tilde \eta$ is purely imaginary with $\tilde\eta>0$.
Let $\hat{x} = (\cos \hat{\theta},\sin \hat{\theta})$ and $d=(\cos \theta_d,\sin \theta_d)$ be the observation and incident directions, respectively.
By the Jacobi-Anger expansion(see e.g., \cite[Formula (3.89)]{kress.1998}), we have
 \begin{equation}
 	e^{i k x \cdot d}=\sum_{n=-\infty}^{\infty} i^{n} J_{n}(k|x|) e^{i n \theta}, \quad x \in \mathbb{R}^{2}.
 \end{equation}
where $\theta=\hat{\theta}-d$ denotes the angle between $\hat{x}$ and $d$.
Using the impedance boundary condition \eqref{imp-bound}, we can get the scattered field $u^s=u^s(x; B_R, d, k, \eta)$ by
\begin{equation}
	u^s(x; B_R, d, k, \eta) = -\sum_{n=-\infty}^{\infty} i^{n}\frac{k J^{\prime}_{n}(kR)+ \eta J_{n}(kR)}{k H^{(1)\prime}_{n}(kR)+ \eta H^{(1)}_{n}(kR)} H^{(1)}_{n}(k |x|) e^{i n \theta}.
\end{equation}
This leads to the far-field pattern $u^{\infty}(\hat{x}; B_R, d, k, \eta)$ of the disk $B_R$:
\begin{equation}
		u^{\infty}(\hat{x}; B_R, d, k, \eta) = -C \sum_{n=-\infty}^{\infty} \frac{k J^{\prime}_{n}(kR)+ \eta J_{n}(kR)}{k H^{(1)\prime}_{n}(kR)+ \eta H^{(1)}_{n}(kR)} e^{i n \theta},
\end{equation}
where $C=\sqrt{\frac{2}{k \pi}} e^{-i \frac{\pi}{4}}$.
Then, an eigen system $(\lambda^{(n)}_{B_R}, \varphi^{(n)}_{B_R})$ of the far-field operator $F_{B_R,imp}$ is given by
  \begin{equation}\label{eig_sys_diskBR_imp}
 \lambda^{(n)}_{B_R}=-2 \pi C \frac{k J^{\prime}_{n}(kR)+ \eta J_{n}(kR)}{k H^{(1)\prime}_{n}(kR)+ \eta H^{(1)}_{n}(kR)}, \qquad  \varphi^{(n)}_{B_R}(\hat{x})=e^{i n \hat\theta}.
 \end{equation}
By the asymptotic behaviour of Bessel functions (see \cite{kress.1998})
 \begin{equation}\label{jn}
 J_{n}(k R) = \frac{(k R)^{n}}{2^{n} n !}{(1+O(\frac{1}{n}))},\quad n \rightarrow +\infty,
 \end{equation}
 \begin{equation}\label{hn}
 H_{n}^{(1)}(k R) = \frac{2^{n}(n-1)!}{\pi i(k R)^{n}}{(1+O(\frac{1}{n}))},\quad n \rightarrow +\infty,
 \end{equation}
 we have
 \begin{equation}\label{eq:4}
 \begin{aligned}
  	|\lambda^{(n)}_{B_R}| &= 2 \pi |C| \left| \frac{k\frac{(k R)^{n-1}}{2^n (n-1)!} + \eta \frac{(k R)^n}{2^n n!}}{k\frac{-2^n n!}{\pi i (k R)^{n+1}} + \eta \frac{2^n (n-1)!}{\pi i (k R)^n}} \right| (1+O(\frac{1}{n}))\\
  	&= 2 \pi \sqrt{\frac{2}{k \pi}} \frac{\pi (k R)^{2n}}{2^{2n} n! (n-1)!} \left| \frac{\eta R + n}{\eta R - n} \right| (1+O(\frac{1}{n}))\\
  	&= \sqrt{\frac{2}{k \pi}} \frac{ \pi^2 (k R)^{2n}}{2^{2n-1} n! (n-1)!}  (1+O(\frac{1}{n}))
 \end{aligned}
 \end{equation}
 and
 \begin{equation}\label{int-circle}
 	\left|\left\langle e^{-i k \hat{x}\cdot z}, \varphi^{(n)}_{B_R}(\hat{x})\right\rangle_{\mathbb{S}}\right|^{2} = \left| J_{n}(k |z|) \right|^2
 	= \frac{(k |z|)^{2n}}{2^{2n} n! n!}(1+O(\frac{1}{n})).
 \end{equation}
 Then we get
 \begin{equation}\label{eq:7}
 	I_{imp}(z) = \left[ \sum_{j} \sqrt{\frac{k}{8 \pi^3}} \frac{1}{j} \left( \frac{|z|}{R} \right)^{2j} (1+O(\frac{1}{j}))\right]^{-1}.
 \end{equation}
  Obviously, the indicator function $I_{imp}(z)>0$ if and only if $|z|<R$, i.e., $z$ lies inside $B_R$. Moreover, recalling the series
 expansion $-\ln \left(1-t\right)= \sum^{+\infty}_{n=1}\frac{t^n}{n}$ for $|t|<1$, the principle part of $I_{imp}(z)$ for $|z|<R$ can be written as
\begin{equation}
\begin{aligned}
	\sum_{j} \sqrt{\frac{k}{8 \pi^3}} \frac{ \left( \frac{|z|^2}{R^2} \right)^{j}}{j} &= -\sqrt{\frac{k}{8 \pi^3}} \ln \left(1-\frac{|z|^2}{R^2}\right) \\
	&= \sqrt{\frac{k}{8 \pi^3}} \ln \left( \frac{R^2}{R^2-|z|^2} \right)\\
	&= -\sqrt{\frac{k}{8 \pi^3}} (\ln (R-|z|) + \ln 2 -\ln R).
\end{aligned}
\end{equation}
This implies that $I^{-1}_{imp}(z)\sim -\ln (R-|z|)$ as $|z|\rightarrow R$, $|z|<R$.

If $B_R$ is a sound-soft disk, we have
 \begin{equation}
 u^\infty(\hat{x}; B_R, d, k) = -C \sum\limits^\infty_{n=-\infty}\frac{J_{n}(kR)}{H^{(1)}_{n}(kR)} e^{in(\hat{\theta}-\theta_d)}.
 \end{equation}
and the spectrum system
 \begin{equation}\label{eig_sys_diskBR}
 \lambda^{(n)}_{B_R}=-2 \pi C \frac{J_{n}(k R)}{H_{n}^{(1)}(kR)}, \quad
 \varphi^{(n)}_{B_R}(\hat{x})=e^{i n \hat\theta}.
 \end{equation}
 From \eqref{jn} and \eqref{hn}, we obtain
 \begin{equation}\label{eq:5}
\begin{aligned}
|\lambda^{(n)}_{B_R}| &= 2 \pi |C| \left| \frac{\frac{(k R)^n}{2^n n!}}{\frac{2^n (n-1)!}{\pi i (k R)^n}} \right| (1+O(\frac{1}{n}))\\
&= 2 \pi \sqrt{\frac{2}{k \pi}} \frac{\pi (k R)^{2n}}{2^{2n} n! (n-1)!} (1+O(\frac{1}{n}))\\
&= \sqrt{\frac{2}{k \pi}} \frac{ \pi^2 (k R)^{2n}}{2^{2n-1} n! (n-1)!}  (1+O(\frac{1}{n})).
\end{aligned}
\end{equation}
This together with \eqref{int-circle} gives
 \begin{equation}
I(z) = \left[ \sum_{j} \sqrt{\frac{k}{8 \pi^3}} \frac{1}{j} \left( \frac{|z|}{R} \right)^{2j} (1+O(\frac{1}{j}))\right]^{-1},
\end{equation}
which is the same as (\ref{eq:7}).
Hence, we can again conclude that $I(z)>0$ if and only if $|z|<R$ and
$I^{-1}(z)\sim -\ln (R-|z|)$ as $|z|\rightarrow R$, $|z|<R$.
\begin{remark}
In three dimensions, it holds that $I(z)\sim (R-|z|)^{-1}$ as $|z|\rightarrow R$, $|z|<R$; see \cite[Chapter 1.5]{Kirsch.2008}. For general sound-soft obstacles,
it holds that $I^{-1}(z)\sim ||\Phi(\cdot, z)||^2_{H^{1/2}(\partial D)}$ as $z\rightarrow\partial D$, $z\in D$.
\end{remark}

\section{Factorization method with one plane wave: a model-driven and data-driven approach}\label{sec:3}

\subsection{Further discussions on Kirch's Factorization method}\label{sub3.1}
Before stating the one-wave version of the factorization method, we first present a corollary of Theorems \ref{class_fac} and \ref{class_fac_imp}. Denote by $\Omega\subset \R^2$ a convex Lipschitz domain which may represent a sound-soft or impedance scatterer. From numerical point of view, $\Omega$ will play the role of test domains for imaging the unknown scatterer $D$. Here we use the notation $\Omega$ in order to distinguish from our target scatterer $D$.
The far-field operator corresponding to $\Omega$ is therefore given by
\begin{equation}
	(F_{\Omega} g)(\hat{x}) = \int_{\mathbb{S}} u^{\infty}_{\Omega} (\hat{x}, d) g(d) ds(d), \quad F_{\Omega} : L^2(\mathbb{S}) \rightarrow  L^2(\mathbb{S}),
\end{equation}
where $u^{\infty}_{\Omega} (\hat{x}, d) $ is the far-field pattern corresponding to the plane wave $e^{ikx\cdot d}$ incident onto $\Omega$. The eigenvalues and eigenfunctions of $F_{\Omega}$ will be denoted by $(\lambda^{(j)}_{\Omega}, \varphi^{(j)}_{\Omega})$.

\begin{corollary}\label{fac_convex_domain} Let $v^\infty\in L^2(\s)$ and assume that
 $k^2$ is not an eigenvalue of $-\Delta$ over $\Omega$ with respect to the boundary condition under consideration. Then 
\begin{equation}
	I(\Omega) = \sum_{j} \frac{\left|\left\langle v^{\infty}, \varphi^{(j)}_{\Omega}\right\rangle_{\mathbb{S}}\right|^{2}}{\left|\lambda^{(j)}_{\Omega}\right|} < + \infty
\end{equation}
	if and only if $v^{\infty}$ is the far-field pattern of some radiating solution $v^{s}$, where $v^{s}$ satisfies the Helmholtz equation
\begin{equation}\label{helmholtz}
	\bigtriangleup v^{s} + k^{2} v^{s} = 0 \qquad \mbox{in} \quad  \mathbb{R}^2 \backslash \overline{\Omega},
\end{equation}
 with the boundary data $v^{s} |_{\partial \Omega}\in H^{1/2}(\partial \Omega).$
\end{corollary}
\begin{remark}
If $v^\infty(\hat x)=e^{-ik \hat x\cdot z}$, then by Theorems \ref{class_fac} and \ref{class_fac_imp} it holds that $I(\Omega)<\infty$ if and only if $z\in \Omega$. This implies that the scattered field $v^s(x)=\Phi(x,y)$ is a well-defined analytic function in $\R^2\backslash\overline{\Omega}$ and
$v^{s} |_{\partial \Omega}\in H^{1/2}(\partial \Omega)$ if and only if $I(\Omega)<\infty$.
 Hence, the results in Corollary  \ref{fac_convex_domain} follow directly from Theorems \ref{class_fac} and \ref{class_fac_imp} in this special case.
\end{remark}

{\bf{Proof}.} In $\Omega$ is sound-soft, by \eqref{fac-D} we have
$F_{\Omega} = -G_{\Omega} S^*_{\Omega} G^*_{\Omega}$,
where $G_{\Omega}: H^{1/2}(\partial \Omega)\to L^2(\mathbb{S})$ is the data-to-pattern operator corresponding to $\Omega$.
Obviously, $I(\Omega) <+\infty$ if and only if $v^{\infty} \in Range((F^*_{\Omega} F_{\Omega})^{1/4})$. Since $Range((F^*_{\Omega} F_{\Omega})^{1/4}) = Range(G_{\Omega})$, we get $v^{\infty} \in Range(G_{\Omega})$ if and only if $I(\Omega) <+\infty$.
Recalling the definition of $G_{\Omega}$, it follows that $v^{s}$ satisfies the Helmholtz equation \eqref{helmholtz} and the Sommerfeld radiation condition \eqref{SommerfeldCond} with the boundary data $v^{s} |_{\partial \Omega}\in H^{1/2}(\partial \Omega).$
The impedance case can be proved in an analogous manner by applying the factorization
$F_{\Omega} = -G_{\Omega,imp} T^*_{\Omega,imp} G^*_{\Omega,imp}$ (see (\ref{fac-Dimp})) and the range identity $Range(F_{\Omega,\#}) = Range(G_{\Omega})$.
$\hfill\square$

\subsection{Explicit examples when $\Omega$ is a disk}
 Corollary \ref{fac_convex_domain} relies essentially on the factorization form (see e.g., (\ref{fac-D}),(\ref{fac-Dimp})) of the far-field operator and is applicable to both penetrable and impenetrable scatterers $\Omega$. However, establishing the abstract framework of the factorization method turns out to be nontrivial in some cases, for instance, time-harmonic acoustic scattering from mixed obstacles and electromagnetic scattering from perfectly conducting obstacles.  Below we show that the results of Corollary \ref{fac_convex_domain} can be justified independently of the factorization form, if the test domain $\Omega$ is chosen to be a disk of acoustically Dirichlet or impedance type. This is mainly  due to the explicit form of far-field patterns for Dirichlet and impedance disks.

Let $B_{h}(z)$ be the disk centered at $z\in \R^2$ with radius $h>0$. The boundary of $B_h(z)$ is denoted by  $\Gamma_{z,h} := \{x: |x-z|=h\}$. It is supposed that $B_h(z)$ is either a sound-soft disk in which $k^2$ is not the Dirichlet eigenvalue of $-\Delta$,  or an impedance disk with the constant impedance coefficient $\eta\in \C$ such that ${\rm Im}(\eta)>0$.
Denote by $u^{\infty}_{z,h}(\hat{x}, d)=u^\infty(\hat{x}; B_h(z), d, k)$ the far-field pattern incited by the plane wave $e^{ikx\cdot d}$ incident onto $B_{h}(z)$ and by $F_{z,h}$ the associated far-field operator, that is,


	\begin{equation}
	(F_{z,h} g)(\hat{x})=\int_\mathbb{S} u^{\infty}_{z,h}(\hat{x}, d) g(d) d s(d) \quad \text { for } \hat{x} \in \mathbb{S},\quad g\in L^2(\s).
	\end{equation}
Using the translation formula
	\begin{equation}
	u^\infty(\hat{x}; B_h(z), d, k) =	e^{ikz\cdot(d-\hat{x})} u^\infty(\hat{x}; B_h(O), d, k),
	\end{equation}
together with the spectral system \eqref{eig_sys_diskBR} and \eqref{eig_sys_diskBR_imp} for $B_h(O)$,
we can get the spectral system $(\lambda^{(n)}_{z,h}, \varphi^{(n)}_{z,h})$ of $F_{z,h}$ under the Dirichlet or impedance boundary condition:
\begin{itemize}
	\item If $B_h(z)$ is a sound-soft disk, then
\begin{equation}\label{eigsys-sf}
{\lambda^{(n)}_{z,h}=-2 \pi C \frac{J_{n}(k h)}{H_{n}^{(1)}(k h)}}, \quad
{\varphi^{(n)}_{z,h}(\hat x)=e^{i n \hat\theta-i k z\cdot(\cos\hat\theta, \sin\hat\theta)}}.
\end{equation}	
In particular, $\lambda^{(n)}_{z,h}\neq 0$ if $k^2$ is not the Dirichlet eigenvalue of $-\Delta$ in $B_h(z)$.
	\item If $B_h(z)$ is an impedance disk with the impedance constant $\eta\in \C$, then
\begin{equation}\label{eigsys}
{\lambda^{(n)}_{z,h}=- 2 \pi C \frac{k J_{n}^{\prime}(k h) + \eta J_{n}(k h)}{k H_{n}^{(1) \prime}(k h) + \eta H_{n}^{(1)}(k h)}}, \quad
{\varphi^{(n)}_{z,h}(\hat x)=e^{i n \hat\theta-i k z\cdot(\cos\hat\theta, \sin\hat\theta)}}.
\end{equation}
In particular, we have $\lambda^{(n)}_{z,h}\neq 0$.
\end{itemize}



 Note that the above eigenvalues are independent of $z$ and the eigenfunctions are independent of $h$.
 Taking $\Omega=B_h(z)$, we can rewrite Corollary \ref{fac_convex_domain} as
		\begin{corollary}\label{thm_Izh}
Let $v^\infty\in L^2(\s)$. Then
		\begin{equation}
		I(z,h) := \sum_{j} \frac{\left|\left\langle v^{\infty}, \varphi^{(j)}_{z,h}\right\rangle_{\mathbb{S}}\right|^{2}}{\left|\lambda^{(j)}_{z,h}\right|} < + \infty
		\end{equation}
		if and only if $v^{\infty}$ is the far-field pattern of the radiating solution $v^{s}$, where $v^{s}$ satisfies Helmholtz equation
		$\Delta v^{s} + k^{2} v^{s} = 0$ in  $|x-z|>h$ with
 the boundary data $f:=v^{s} |_{\Gamma_{z,h}}\in H^{1/2}(\Gamma_{z,h})$.
	\end{corollary}

    	{\bf{Proof}.} Without loss of generality, we assume that the center $z$ coincides with the origin, i.e., $B_{h}(z)=B_{h}(O)$. By the Jacobi-Anger expansion (see e.g.,\cite{kress.1998}), $v^s$ can be expanded into the series
    	\begin{equation} \label{expend_scatterfieid}
    	v^{s}(x) = \sum_{n \in \mathbb{Z}} A_{n} H_{n}^{(1)}(k |x|) e^{i n \hat\theta},\quad |x|>R,
    \quad x=(|x|,\hat\theta),
    	 \end{equation}
    	 for some sufficiently large $R>0$, with the far-field pattern
given by (see \cite[(3.82)]{kress.1998})
    \begin{equation}
    v^{\infty}(\hat{x}) = \sum_{n \in \mathbb{Z}} A_{n} C_{n} e^{in \hat\theta}, \quad C_{n}:=\sqrt{\frac{2}{k \pi}} e^{-i (\frac{n\pi}{2}+\frac{\pi}{4})}.
    \end{equation}
Recall from the asymptotic behavior \eqref{eq:4} and \eqref{eq:5} that
    \be\label{lambda}
      \left|\lambda^{(j)}_{z,h}\right| 
    =2 \pi ^2 \sqrt{\frac{2}{\pi k}} \frac{(k h)^{2j}}{2^{2j} j! (j-1)!}{(1+O(\frac{1}{j}))}
    \en
as $j\rightarrow\infty$.
 Hence, the function $I(z,h)$ with $z=O$ can be written as
    \begin{equation}\label{Izh}
    \begin{split}
    	I(z,h) &=\sum_{j \in \mathbb{Z}} \frac{\left|\left\langle  \sum A_{n} C_{n} e^{in \hat\theta} , \varphi^{(j)}_{z,h}\right\rangle_{\mathbb{S}}\right|^{2}}{\left|\lambda^{(j)}_{z,h}\right|} \\
    	&=\sum_{j \in \mathbb{Z}} \left|\lambda^{(j)}_{z,h}\right|^{-1} |2\pi A_j C_j|^2 \\
    	&= \sum_{j \in \mathbb{Z}} \sqrt{\frac{8}{\pi k}} |A_{j}|^{2} \frac{2^{2j} j! (j-1)!}{(k h)^{2j}} {(1+O(\frac{1}{j}))}.
   	\end{split}
    \end{equation}
    On the other hand, by \eqref{expend_scatterfieid} we have
    \begin{equation}\label{vs}
    v^s(x) =\sum_{j \in \mathbb{Z}} |A_{j}|\frac{2^{j}(j-1)!}{\pi (k|x|)^{j}}{(1+O(\frac{1}{j}))},     \end{equation}
    By the definition of the norm $ ||\cdot||_{H^{1/2}(\Gamma_{z,h})} $, it is easy to see when $z=O$  that
\begin{equation}\label{vs_bound}
    \begin{split}
       ||v^{s}||^{2}_{H^{1/2}(\Gamma_{z,h})} &= \sum_{j \in \mathbb{Z}} (1+j^2)^{1/2} |A_j H^{(1)}_j(kh)|^2\\
       &=\sum_{j \in \mathbb{Z}} |A_{j}|^{2} \frac{2^{2j} j! (j-1)!}{\pi ^2 (k h)^{2j}}{(1+O(\frac{1}{j}))}.
    \end{split}
    \end{equation}
Obviously, the series \eqref{Izh} and \eqref{vs_bound} have the same convergence. On the other hand, by \cite[Theorem 2.15]{kress.1998},
the boundedness of $||v^s||_{L^2(\Gamma_{z,h})}$ means that $v^s$ is a radiating solution in $|x-z|>h$ with the far-field pattern $v^\infty$. This proves that $I(z,h)<\infty$ if and only if $v^s$ is a radiating solution in $|x-z|>h$ with the far-field pattern $v^\infty$ and with the $H^{1/2}$-boundary data on $\Gamma_{z,h}$. $\hfill\square$


In our applications of Corollaries \ref{fac_convex_domain} and \ref{thm_Izh}, we will take $v^\infty$ to be the measurement data $u_D^\infty(\hat x; d_0):=u^\infty(\hat x; d_0, D)$ that corresponds to our target scatterer $D$ and the incident plane wave $e^{ikx\cdot d_0}$ for some fixed $d_0\in\s$. We shall omit the dependance on $d_0$ if it is always clear from the context.  Our purpose is to extract the geometrical information on $D$ from the domain-defined indicator function $I(\Omega)$ or $I(z,h)$. By Corollary \ref{fac_convex_domain}, $I(\Omega)<\infty$ if the scattered field $u^s_D(x)=u^s(x; d_0, D)$ can be extended to the domain $\R^2\backslash\overline{\Omega}$ as a solution to the Helmholtz equation. This implies that $u^s_D$ admits an analytical extension across the boundary $\partial D$ when the inclusion relation $D\subset\Omega$ does not hold. Therefore,
the one-wave factorization method requires us to exclude the possibility of analytical extension, which however is possible only if $\partial D$ is not everywhere analytic. If $\partial D$ is analytic, it follows from the the Cauchy-Kovalevski theorem (see e.g. \cite[Chapter 3.3]{John}) that $u_D^s$ can be locally extended into the inside of $D$ across the boundary $\partial D$. In the special case that $D=B_h(O)$, by the Schwartz reflection principle the scattered field $u_D^s$ can be globally continued into $\R^2\backslash\{O\}$.
Below we shall discuss the absence of the analytic extension of $u^s_D$ in corner domains, which is an import ingredient in establishing the one-wave factorization method.

\subsection{Absence of analytical extension in corner domains}

We first consider time-harmonic acoustic scattering from a convex polygon of sound-soft, sound-hard  or impedance type. In the impedance case, the impedance function is supposed to be a constant.

\begin{lemma}\label{lemma_obstacle}
	Assume that $D$ is either a sound-soft, sound-hard or impedance obstacle occupying a convex polygon. Then the scattered field $u^s_D(x,d_0)$ for a fixed $d_0\in \s$ cannot be analytically extended from $\mathbb{R}^{2} \backslash \overline{D}$ into $D$ across any corner of $D$.
\end{lemma}

As one can imagine, the proof of Lemma \ref{lemma_obstacle} is closely related to uniqueness in determining a convex polygonal obstacle with a single incoming wave (see e.g., \cite{cheng2004},\cite[Theorem 5.5]{kress.1998} and \cite{HM}). In fact, the result of Lemma \ref{lemma_obstacle} implies that a convex polygonal obstacle of sound-soft, sound-hard or impedance type can be uniquely determined by one far-field pattern. There are several approaches to prove Lemma \ref{lemma_obstacle}. Below we present a unified method valid for any kind boundary condition under consideration. The original version of this approach was presented in \cite{FI89} for proving uniqueness in inverse conductivity problems.

{\bf Proof.} Assume on the contrary that $u^s_D$ can be analytically continued across a corner of $\partial D$. By coordinate translation and rotation, we may suppose that this corner coincides with the origin, so that $u^s_D$ and also the total field $u_D=u^s_D+u^i$ satisfy the Helmholtz equation in $B_\epsilon(O)$ for some $\epsilon>0$. Since $u_D$ is real analytic in $(\R^2\backslash\overline{D})\cup B_\epsilon(O)$ and $D$ is a convex polygon, $u_D$ satisfies the Helmholtz equation in a neighborhood of an infinite sector $\Sigma\subset \R^2\backslash\overline{D}$ which extends the finite sector $B_\epsilon(O)\cap D$ to $\R^2\backslash\overline{D}$. On the other hand, $u_D$ fulfills the boundary condition on the two half lines $\partial \Sigma$ starting from the corner point $O$. 
By the Schwartz reflection principle of the Helmholtz equation, $u_D$ can be extended onto $\R^2$ (see also \cite{FI89} for discussions on the conductivity equation). We remark that $u_D\equiv 0$ when the angle of $\Sigma$ is irrational  and that the impedance case follows from the arguments in \cite{HM}. This implies that the scattered field $u_D^s$ is an entire radiating solution. Hence, $u_D^s$ must vanish identically in $\R^2$ and $u_D=u^i$ must fulfill the boundary condition on $\partial D$. However, this is impossible for a plane wave incidence under either of the Dirichlet,Neumann or Robin condition.
$\hfill\square$

Next we consider the source radiating problem where $D$ is a polygonal source term. The absence of analytical extension in this case implies that a polygonal source term cannot be a non-radiating source (that is, the resulting far-field pattern vanishes identically). This can be proved based on the idea of constructing CGO solutions \cite{BLS,Bl2018,HSV} or analyzing corner singularities for solutions of the inhomogeneous Laplace equation \cite{hu2018,HL}. Below we consider a  special analytic source function, which will significantly simplify the arguments employed in \cite{ElHu2015,hu2018} for inverse medium scattering problems and those in
\cite{Bl2018,HL} for inverse source problems.
\begin{lemma}\label{lemma_source_const}
Assume that $D$ is a convex polygon and let $\chi_D$ be the characteristic function for $D$. If $u\in H^2_{loc}(\R^2)$ is a radiating solution to
\be\label{source}
\Delta u(x)+k^2u(x)=\chi_D(x)f(x)\quad\mbox{in}\quad \R^2,
\en
where $f(x)$ is real-analytic and non-vanishing near the corner point $O$ of $D$ and the lowest order Taylor expansion of $f$ at $O$ is harmonic.
Then $u$ cannot be analytically extended from $\R^2\backslash\overline{D}$ to $D$ across the  corner $O$. 
\end{lemma}

\begin{figure}[ht]
	\centering
	\includegraphics[scale=0.3]{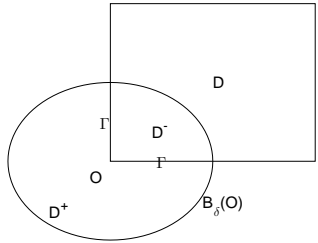}   
	\caption{Illustration of a convex polygonal source term where $O$ is corner point of $D$.}\label{pic2}
\end{figure}
{\bf Proof.} Assume that $u$ can be analytically extended from $\R^2\backslash\overline{D}$ to $B_\delta(O)$ for some $\delta > 0$ (see Figure \ref{pic2}), as a solution of the Helmholtz equation. Without loss of generality, the corner point $O$ is supposed to coincide with the origin. Set $v^\pm=u|_{D^\pm}$ where $D^+:=B_\delta(O)\cap (\R^2\backslash\overline{D})$ and $D^-:=B_\delta(O)\cap D$. This implies that
$\Delta v^+ + k^2v^+ = 0$ in $B_\delta(O)$ and thus the Cauchy data $v^-=v^+$, $\partial_\nu v^-=\partial_\nu v^+$ on $\Gamma:=\partial D\cap B_\delta(O)$ are analytic. Since $f$ is also analytic in $B_\delta(O)$,
by the Cauchy-Kovalevskaya theorem (see e.g. \cite[Chapter 3.3]{John}) the function $v^-$ can also be analytically continued into $B_\delta(O)$ as a solution of $\Delta v^- + k^2v^- = f$ in $B_\delta(O)$.
Setting $w=v^--v^+$ in $B_\delta(O)$, we have
\be\label{lr}
		\Delta w + k^2w = f\quad\mbox{in}\quad B_\delta(O),\quad
		w=\partial_{\nu}w=0\quad\mbox{on}\quad \Gamma.
	\en
Using \cite[Lemma 2.2]{EHY2015}, the analytic functions $w$ and $f$ can be expanded in polar coordinates into the series
\be
\begin{split}\label{expansion}
	&w(x)=\sum_{n+2m \geq 0}r^{n+2m}(a_{n,m}\cos n\theta+b_{n,m}\sin n\theta),\\
&f(x)=r^{N}(\tilde a_{N}\cos N\theta+\tilde b_{N}\sin N\theta)+\sum_{n+2m \geq N+1}r^{n+2m}(\tilde a_{n,m}\cos n\theta+\tilde b_{n,m}\sin n\theta),
\end{split}\en
in $B_\delta(O)$,
where $(r,\theta)$ denote the polar coordinates of $x\in \R^2$ and $N\in \N_0$. Since $f$ does not vanish identically, we may suppose that $|\tilde a_{N}|+|\tilde b_{N}|>0$.
Recalling the Laplace operator in polar coordinates, $\Delta = \frac{\partial^2}{\partial r^2}+\frac{1}{r}\frac{\partial}{\partial r}+\frac{1}{r^2}\frac{\partial^2}{\partial \theta^2}$, we have
\ben
f(x)&=&\Delta w(x)+k^2w(x) \\
		&=&\sum_{n+2m \geq 0}[4(n+m+1)(m+1)r^{n+2m}(a_{n,m+1}\cos n\theta+b_{n,m+1}\sin n\theta)]\\
		&&+\sum_{n+2m \geq 0} k^2r^{n+2m}(a_{n,m}\cos n\theta+b_{n,m}\sin n\theta) \\
\enn
Inserting the expansion of $f(x)$ in \eqref{expansion} and comparing the coefficients of $r^{l}$, $l\in\N_0$, we can get the recurrence relations for $a_{n,m}$,
\begin{equation}\label{recursion_cx}
	\begin{aligned}
		4(N+1) a_{N,1} +k^2 a_{N,0} &= \tilde{a}_{N}, \\
		4(n+m+1) (m+1)a_{n,m+1} +k^2 a_{n,m} &= 0, \quad n+2m <N,\\
		4(n+m+1) (m+1)a_{n,m+1} +k^2 a_{n,m} &= 0, \quad n+2m = N, m>0.
	\end{aligned}
\end{equation}
The same relations hold for $b_{n,m}$.
Now, we suppose without loss of generality that $\Gamma=\{(r,\pm\theta_0): |r|<\delta\}$ for some $\theta_0\in(0,\pi/2)$.
From the boundary conditions $w=\partial_{\theta}w=0$ on $\Gamma$, we have
\begin{equation}\label{boundformula}
	\left\{
	\begin{array}{lr}
		\sum\limits_{n,m\in \mathbb{N}_0, n+2m=l}a_{n,m}\cos n\theta_{0} = 0,\\
		\sum\limits_{n,m\in \mathbb{N}_0,n+2m=l}na_{n,m}\sin n\theta_{0} = 0,\\
		\sum\limits_{n,m\in \mathbb{N}_0,n+2m=l}nb_{n,m}\cos n\theta_{0} = 0,\\
		\sum\limits_{n,m\in \mathbb{N}_0,n+2m=l}b_{n,m}\sin n\theta_{0} = 0,
	\end{array}
	\right.
\end{equation}
for any $l \in \mathbb{N}_0$.
%
From the second formula in \eqref{recursion_cx} and the first two formulas in \eqref{boundformula}, we can easily obtain $a_{n,m}=0$ if $n+2m<N$.

Now we prove that $a_{n,m}=0$ if $l=n+2m=N$.
In fact, for $m\geq 1$, setting $m'=m-1\geq 0$ we derive from the second formula in \eqref{recursion_cx}
that
\ben
a_{n,m}=a_{n,m'+1}=-\frac{k^2}{4(n+m)m}a_{n,m'}=0,
\enn since $n+2m'<N$. The above relations together with the first two formulas in (\ref{boundformula}) with $l=N$ lead to
\ben
a_{N,0}\cos N\theta_{0} = 0,\quad
		Na_{N,0}\sin N\theta_{0} = 0,
\enn which imply that $a_{N,0}=0$.

%
%
When $l =n+2m= N+2$, we observe that $n+2m'=N$, where $m'=m-1$. Hence one can get $a_{n,m}=0$ if $m>1$ and $l=N+2$, by using the third formula in \eqref{recursion_cx} and the fact that $a_{n,m'}=0$ for all $n+2m'=N$. Then it follows from the first two formulas in (\ref{boundformula}) with $l=N+2$ that
\begin{equation}
	\left\{\begin{array}{lll}
		a_{N+2,0}\cos (N+2)\theta_{0} + a_{N,1}\cos N\theta_{0}= 0,\\
		(N+2)a_{N+2,0}\sin (N+2)\theta_{0} +Na_{N,1}\sin N\theta_{0} = 0,
	\end{array}\right.
\end{equation}
Since $0<\theta_0<\frac{\pi}{2}$, we have (see \cite{ElHu2015})
\begin{equation}
	\left|
	\begin{aligned}
		&\cos (N+2)\theta_{0} \quad &\cos N\theta_{0}\\
		&(N+2)\sin (N+2)\theta_{0} \quad &N\sin N\theta_{0}
	\end{aligned}
	\right|
	= (N+1)\sin 2\theta_{0} -\sin (N+1)\theta_{0} \neq 0.
\end{equation}
Then we have $a_{N+2,0}=a_{N,1}=0$. By the first relation in \eqref{recursion_cx} we get
 $\tilde{a}_{N}=0$. Analogously one can prove $\tilde{b}_{N}=0$. This implies that $f\equiv 0$, which is a contradiction.
$\hfill\square$
\begin{remark}
\begin{itemize}
\item[(i)] If $f$ satisfies the elliptic equation
\ben
\Delta f(x) + A(x) \cdot \nabla f(x)+ b(x)\, f(x) = 0
\enn
where $A(x)=(a_1(x), a_2(x))$ and $b(x)$ are both real-analytic,  then the lowest Taylor expansion of $f$ at any point must be harmonic (see \cite{HL}). Therefore, the class of source functions $f$ specified in  Lemma \ref{lemma_source_const} covers at least harmonic functions, including constant functions.
\item[(ii)] The analyticity of $f$ in Lemma \ref{lemma_source_const} can be weakened to be H\"older continuous near $O$ with the asymptotic behavior (see \cite{HL})
    \begin{equation} \label{eq:f}
		f(x) = r^{N}(A_n \cos N \theta + B_N \sin N\theta) +  o(r^{N}), \quad |x| \to 0,\,
		\end{equation}
		for some $N\in \N_0$ and $A_N, B_N\in\mathbb C$ with $|A_N|+|B_N|>0$. Moreover, the corner  $O\in \partial D$ can be weakened to be a weakly singular point of arbitrary order such that $\partial D$ is not of $C^\infty$-smooth at $O$; see \cite{LHY}.
\end{itemize}
\end{remark}

\subsection{One-wave version of the factorization method}
\begin{figure}[ht]
	\centering
\includegraphics[width=0.8\linewidth]{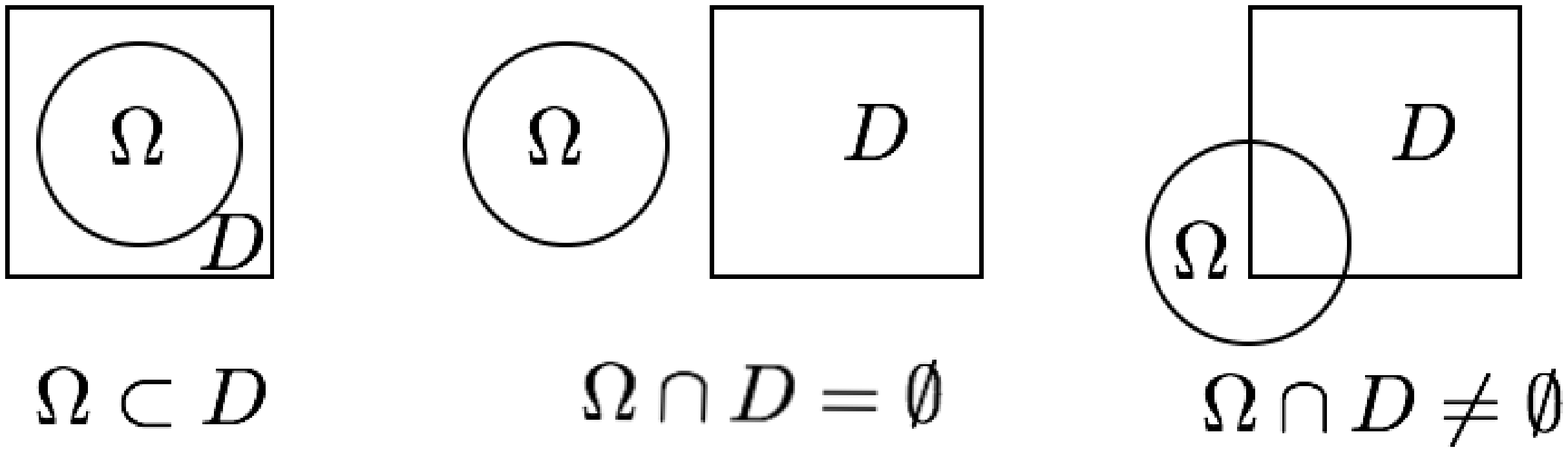}
	\caption{Three cases when $D\subseteq\Omega$ does not hold. Here the test domain $\Omega$ is chosen as a disk in 2D.}
	\label{3case}
\end{figure}

To state the one-wave factorization method, we shall restrict our discussions to convex polygonal impenetrable scatterers of sound-soft, sound-hard or impedance type and to convex polygonal source terms where the source function satisfies the condition of Lemma \ref{lemma_source_const}. In the former case, $u_D^\infty$ represents the far-field pattern of the scattered field caused by some plane wave incident onto $D$; in the latter case, $u_D^\infty$ denotes the far-field pattern of the
radiating solution to \eqref{source}.
Recall from Subsection \ref{sub3.1} that $\Omega$ is a convex sound-soft or impedance scatterer such that $k^2$ is not the eigenvalue of $-\Delta$ in $\Omega$. Denote by $(\lambda^{(j)}_{\Omega}, \varphi^{(j)}_{\Omega})$ the eigenvalues and eigenfunctions of the far-field operator $F_{\Omega}$. Below we characterize the inclusion relationship  between our target scatterer $D$ and the test domain $\Omega$ through the interaction of the measurement data $u_D^\infty$ and the spectra of $F_{\Omega}$.
\begin{theorem}\label{One-FM}
	Define	$$
	W(\Omega) := \sum_{j} \frac{\left|\left\langle u_D^{\infty}, \varphi^{(j)}_{\Omega}\right\rangle_{\mathbb{S}}\right|^{2}}{\left|\lambda^{(j)}_{\Omega}\right|}.
	$$
Then $W(\Omega)<\infty$	if and only if $D \subseteq \Omega$.
\end{theorem}
{\bf Proof.} $\Longrightarrow$ : By Corollary \ref{fac_convex_domain}, $W(\Omega) < +\infty$ implies that $u$ is analytic in $\mathbb{R}^2 \backslash \overline{\Omega}$.
If  $D \nsubseteq \Omega$, three cases might happen (see Fig.\ref{3case}): (i) $\Omega\subset D$; (ii) $\Omega\cap D=\emptyset$; (iii) $\Omega\cap D\neq \emptyset$ and $\Omega\cap (\R^2\backslash\overline{D})\neq \emptyset$.
In either of these cases, we observe that
$u$ can be analytically continued from $\R^2\backslash\overline{D}$ to $D$ across a corner of $\partial D$, which however is impossible by  Lemmas \ref{lemma_obstacle} and \ref{lemma_source_const}. This proves the relationship $D \subseteq \Omega$.

$\Longleftarrow$: We only consider the case where $D$ is an impenetrable scatterer. The source problem can be proved analogously. Assume $D \subseteq \Omega$. Then the scattered field $u^{s}$ satisfies the Helmholtz equation $ \Delta u^{s} + k^{2} u^{s} = 0$ in  $\mathbb{R}^2 \backslash \overline{\Omega}$ with the boundary data  $u^{s}|_{\partial \Omega} \in H^{1/2}(\partial \Omega)$. Then we can get $W(\Omega) < +\infty$ by applying Corollary \ref{fac_convex_domain}.$\hfill\square$

Taking the test domain $\Omega$ as the disk $B_h(z)$, we can immediately get
\begin{corollary}\label{W}
	Let $(\lambda^{(j)}_{z,h}, \varphi^{(j)}_{z,h})$ be an eigensystem of the far-field operator $F_{z,h}$. 	
Define
	\be\label{Wzh}
	W(z,h) := \sum_{j} \frac{\left|\left\langle u_D^{\infty}, \varphi^{(j)}_{z,h}\right\rangle_{\mathbb{S}}\right|^{2}}{\left|\lambda^{(j)}_{z,h}\right|}.
	\en
Then we have $W(z,h)<\infty$
if $h\geq \max_{y\in \partial D}|y-z|$ and $W(z,h)=\infty$
	if $h< \max_{y\in \partial D}|y-z|$.
\end{corollary}

\begin{remark}
Theorem \ref{One-FM}  and Corollary \ref{W} explain how do the a priori data $u_\Omega^\infty(\hat x, d)$ for all $\hat x, d\in \s$ encode the information of the unknown target $D$. Evidently, their proofs rely essentially on mathematical properties of the scattering model. On the other hand, in the terminology of learning theory and data sciences, the test domains $\Omega$ can be regarded as samples and the a priori data $u_\Omega^\infty(\hat x, d)$ are the associated sample data, which are usually calculated off-line. In this sense, the one-wave factorization method is a both model-driven and data-driven approach.
\end{remark}

\begin{figure}[htbp]
		\centering
		\includegraphics[width=0.4\linewidth]{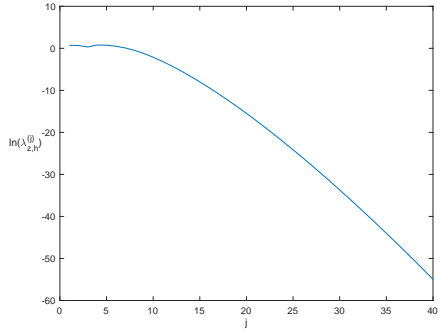}
		\caption{For $k=6$ and $h=1$, $\ln (\lambda^{(j)}_{z,h})$ decays almost linearly as $j\rightarrow\infty$.}
		\label{fig:attenuation of eigenvalues}
	\end{figure}

Corollary \ref{W} says that the maximum distance between a sampling point $z\in \R^2$ and our target $D$ is coded in the function $h\mapsto W(z,h)$. Hence, changing the sampling points $z$ on a large circle $|z|=R$ such that $D\subset B_R(O)$ and computing $\max_{y\in \partial D}|z-y|$ for each $z$ would give an image of $D$ as follows:
\be\label{scheme}
D=\bigcap_{|z|=R,h\in(0,2R)}^{W(h,z)<\infty} B_h(z).
\en
We remark that a proper regularization scheme should be employed in computing the truncated indicator \eqref{Wzh}, because the far-field operator $F_{z,h}$ is compact and the eigenvalues
$\lambda^{(j)}_{z,h}$ decay almost exponentially as $j\rightarrow\infty$; see \eqref{lambda} and Figure \ref{fig:attenuation of eigenvalues}.

\begin{figure}[H]
	\centering    
	\subfigure 
	{
		\begin{minipage}[t]{0.3\linewidth}
			\centering          
			\includegraphics[scale=0.12]{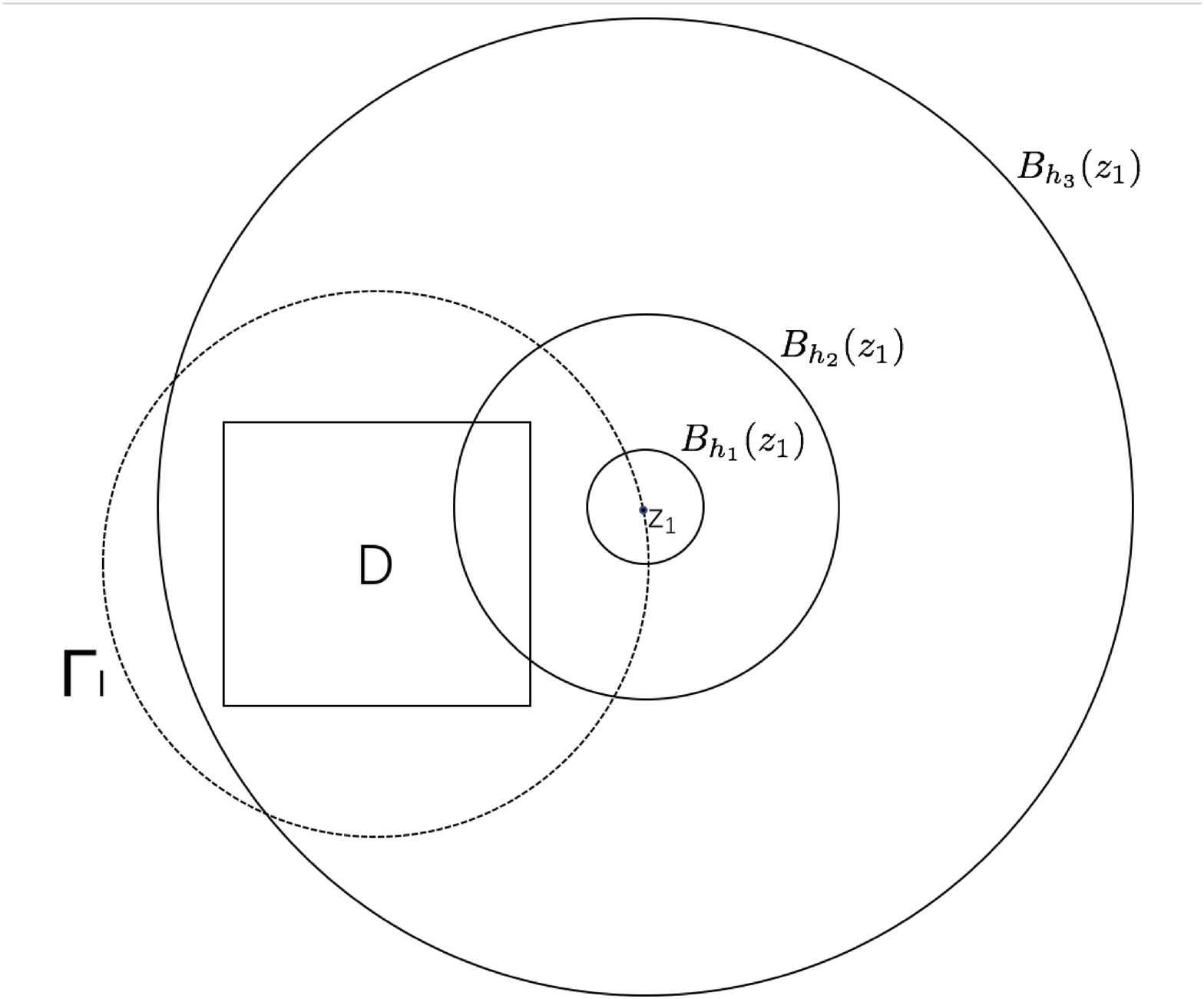}
		\end{minipage}
	}
	\subfigure 
	{
		\begin{minipage}[t]{0.3\linewidth}
			\centering      
			\includegraphics[scale=0.12]{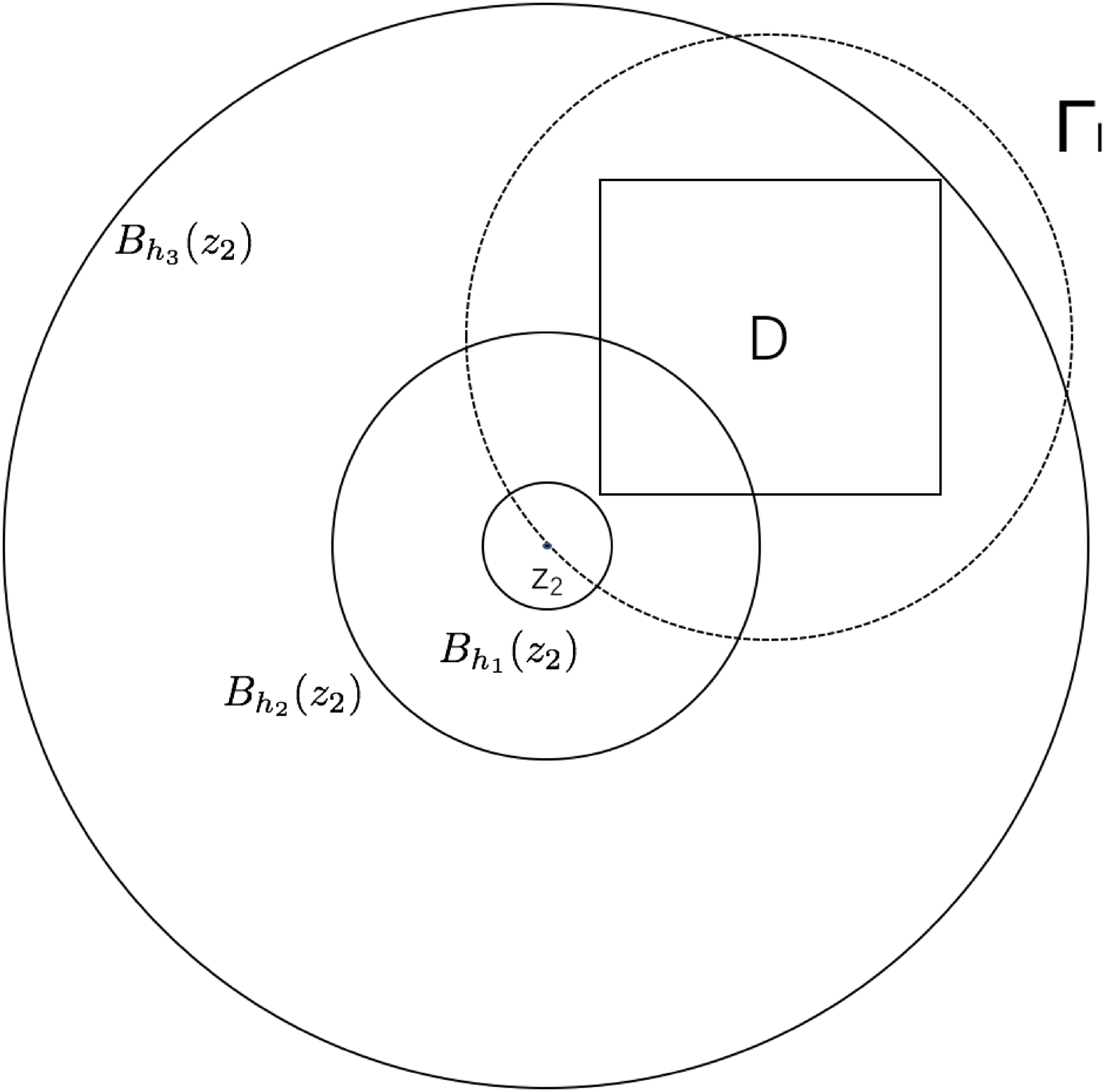}
		\end{minipage}
	}
	\subfigure 
	{
		\begin{minipage}[t]{0.3\linewidth}
			\centering      
			\includegraphics[scale=0.12]{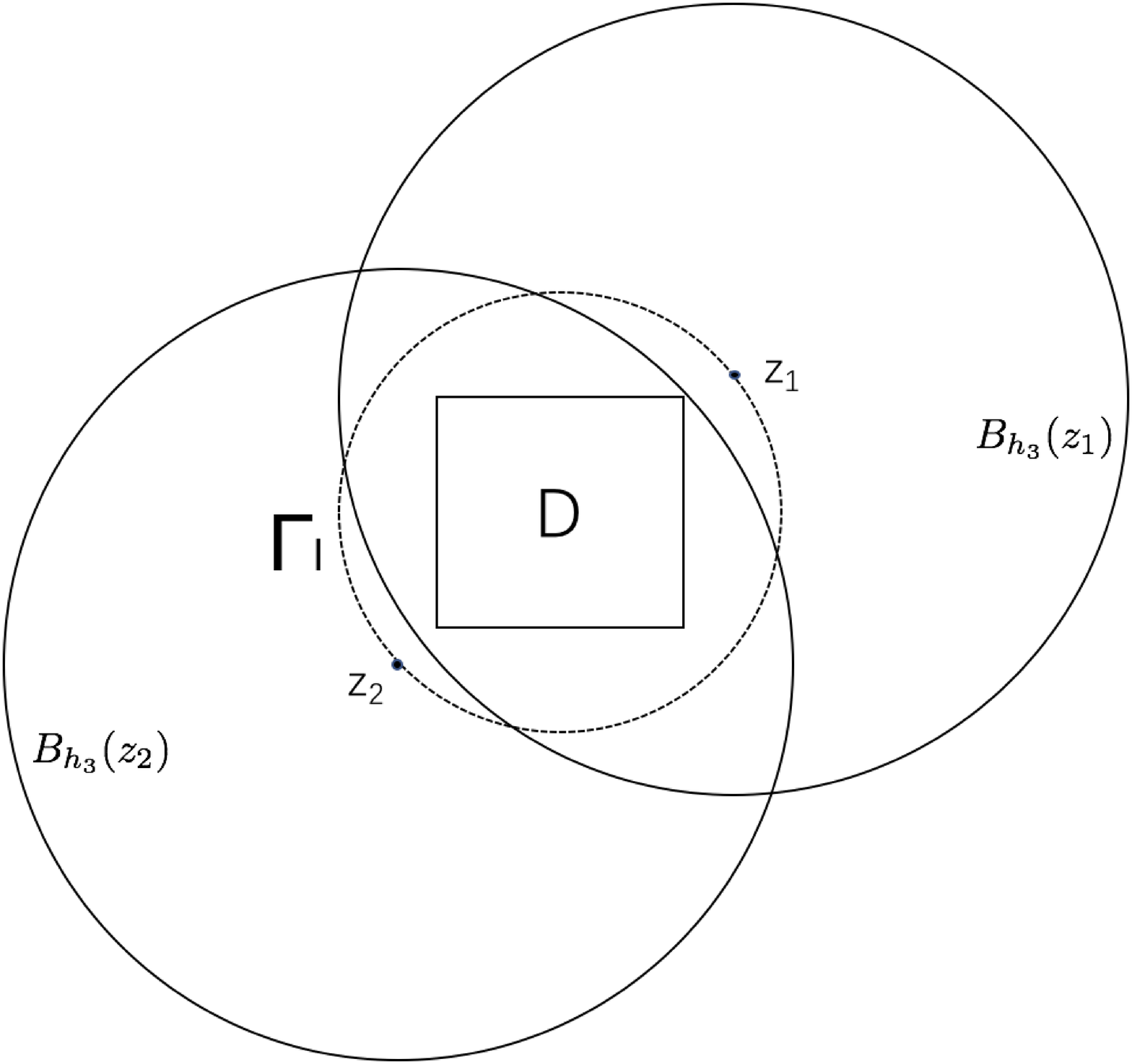}
		\end{minipage}
	}
		\caption{Given $z_1,z_2\in \Gamma_R$, we have $\widetilde{W}(z_m, h_j)=0$ for $m,j=1,2$,
since the underlying obstacle $D$ is not covered by the sampling disks $B_{h_j}(z_m)$ (see the left and middle figures). However, $\widetilde{W}(z_m, h_3)<\infty$ for $m=1,2$, because $D\subset B_{h_3}(z_1)\cap B_{h_3}(z_2)$ (see the right figure).}\label{fig:image}
\end{figure}

If $\Omega=B_h(z)$ is a sound-soft test disk, the $(F^*F)^{1/4}$-method yields the relation $W(z,h)=||g_{z,h}||_{L^2(\s)}^2$, where $g=g_{z,h}\in L^2(\s)$ solves the operator equation

\begin{equation}\label{unTIK}
(F^*_{z,h}F_{z,h})^{1/4}g=u_D^\infty.
\end{equation}
The solution of the equation \eqref{unTIK} is $$g=\sum\limits_{j}\frac{\left\langle u_D^\infty,\varphi_{j}\right\rangle_{\s}}
{\sqrt{\lambda^{(j)}_{z,h}}}\varphi^{(j)}_{z,h}
\qquad\qquad\mbox{if}\quad u_D^\infty\in {\rm Range}((F^*_{z,h}F_{z,h})^{1/4}).$$
Using the Tikhonov regularization we aim to solve the equation
\begin{equation}\label{TIK}
\alpha I+(F^*_{z,h}F_{z,h})^{1/2}g_\alpha=(F^*_{z,h}F_{z,h})^{1/4}u_D^\infty
\end{equation}
with the solution given by \be\label{g}
g_\alpha=\sum\limits_{j}\frac{\sqrt{\lambda^{(j)}_{z,h}}}{|\alpha+\lambda^{(j)}_{z,h}|}\left\langle u_D^\infty,\varphi^{(j)}_{z,h}\right\rangle_{\s} \varphi^{(j)}_{z,h}
\en
where $\alpha>0$ is the regularization parameter.
This implies that
\ben
||g_\alpha||_{L^2(\s)}^2=
\sum\limits_{j}\frac{|\lambda^{(j)}_{z,h}|}
{(|\lambda^{(j)}_{z,h}+\alpha|)^2}\left|\left\langle u_D^\infty,\varphi^{(j)}_{z,h}\right\rangle\right|^2.
\enn
Hence, in our numerics we will use the modified indicator
\be\label{Wt}
\widetilde{W}(z,h) = \left[\sum_{j\leq N} \frac{|\lambda^{(j)}_{z,h}| \left| \left\langle u_D^{\infty}, \varphi^{(j)}_{z,h} \right\rangle_{\mathbb{S}} \right| ^{2}}{\left| \lambda^{(j)}_{z,h}+\alpha \right|^2}\right] ^{-1}
=1/||g_\alpha||_{L^2(\s)}^2.
\en
Our imaging scheme I is described as follows (see Figure \ref{fig:image}):
\begin{itemize}
	\item Suppose that $D\subset B_R(O)$ for some $R>0$ and collect the measurement data $u_D^{\infty}(\hat{x})$ for all $\hat{x} \in \mathbb{S}$;
	\item Choose sampling points $z_{n} \in \Gamma_{R}:= \{x:\ |x|=R\}$ for $n = 1,\cdots,N_{z}$;
	\item Choose $h_m\in(0,2R)$ to get different spectral systems $(\lambda^{(j)}_{z_n,h_m}, \varphi^{(j)}_{z_n,h_m})$ (see (\ref{eigsys-sf}) or (\ref{eigsys}));
	\item For each $z_n\in \Gamma_R$, calculate the maximum distance between $z_n$ and $D$ by $h_{z_n}:=\inf\{h_m\in(0,2R):\widetilde{W}(z_n,h_m)\geq \delta\}$ where $\delta>0$ is a threshhold.
	\item Take $D=\bigcap_{1\leq n\leq N_z}B_{h_{z_n}}(z_{n})$.
\end{itemize}

 In the above numerical scheme, we take the test domain $\Omega$ as sound-soft or impedance disks $B_h(z)$, because the spectral systems $(\lambda^{(j)}_{z,h}, \varphi^{(j)}_{z,h})$ are given explcitly. For a general test domain, the spectral systems $(\lambda^{(j)}_{\Omega}, \varphi^{(j)}_{\Omega})$ should be calculated off-line, so that they are available before inversion. The sample variables (disks) in the above one-wave factorization method consist of the centers $z\in \partial B_R$ and and radii $h\in(0,2R)$. Obviously, the number of these variables is comparable with that of the original factorization method with infinitely many plane waves.

\subsection{Discussions on other domain-defined sampling methods}
There exist some other domain-defined sampling methods in recovering impenetrable scatterers with a single far-field pattern such as range test \cite{rangetest}, no-response test \cite{P2003} and extended sampling method \cite{esm2018}.
It was shown in \cite[Chapter 15]{NP2013} and \cite{LNPW} that range test and no-response test are dual and equivalent for inverse scattering and inverse boundary value problems.
The extended sampling method \cite{esm2018} suggests solving the first kind linear integral equation
\begin{equation}\label{esm}
	(F_{z,h} \tilde g)(\hat{x}) = u^{\infty}_{D} (\hat{x}),
\end{equation}
with regularization schemes. It was proved in \cite{esm2018} that the regularized solution $||\tilde g_\alpha||_{L^2(\s)}^2=\infty$ if $u_D$ cannot be analytically extended into the domain $|x-z|>h$, where $\alpha>0$ is the regularization parameter. We observe that, by
Tikhonov regularization, the solution $\tilde{g}_\alpha=(F_{z,h}^* F_{z,h}+\alpha I)^{-1} F_{z,h}^* u_D^\infty$ to (\ref{esm}) is given by
\ben
\tilde{g}_\alpha=\sum\limits_{j}\frac{\overline{\lambda}^{(j)}_{z,h}}
{|\lambda^{(j)}_{z,h}|^2+\alpha}\left\langle u_D^\infty,\varphi^{(j)}_{z,h}\right\rangle\varphi^{(j)}_{z,h}.
\enn
Obviously, we have
\be\label{esm-i}
||\tilde{g}_\alpha||_{L^2(\s)}^2=
\sum\limits_{j}\frac{|\lambda^{(j)}_{z,h}|^2}
{(|\lambda^{(j)}_{z,h}|^2+\alpha)^2}\left|\left\langle u_D^\infty,\varphi^{(j)}_{z,h}\right\rangle\right|^2
\en
and $||\tilde g_{\alpha} ||_{L^2(\mathbb{S})} < +\infty$ if and only if $D \subset B_{h}(z)$.
Note that in our discussions, $D$ always represents a convex polygonal impenetrable scatterers or a convex polygonal source term. Comparing  \eqref{esm-i} with our regularized solution $g_\alpha$ to \eqref{TIK}, we find for any fixed $\alpha>0$ that,
\begin{equation}
	||\tilde{g}_\alpha||^2_{L^2(\mathbb{S})}  \leqslant  c_\alpha\, ||g_{\alpha}||^2_{L^2(\mathbb{S})},\quad c_\alpha>0,
\end{equation}
since $\lambda^{(j)}_{z,h}\rightarrow 0$ as $j\rightarrow\infty$.

The one-wave version of Range Test considers the first kind integral equation
\begin{equation}\label{range}
	(H^* \varphi)(\hat{x}) = u^{\infty}_{D} (\hat{x}),\quad \varphi\in L^2(\partial \Omega),
\end{equation}
where
$H^*$ is the adjoint operator of the Herglotz operator $H$ with respect to the $L^2$ inner product,  given by
\begin{equation}\label{range-int}
	\left(H^{*} \varphi \right)(\hat{x}):=\int_{\partial \Omega} e^{-i k \hat{x} \cdot y} \varphi (y) d s(y), \quad \hat{x} \in \mathbb{S}.
\end{equation}
Here $\Omega$ is a convex test domain.
Applying Lemmas \ref{lemma_obstacle} and \ref{lemma_source_const}, one can also prove that
 $D \subset \Omega$ if and only if $||\varphi ||_{L^2(\partial \Omega)} < +\infty$. Below we show the equivalence of the domain-defined functionals for the Extended Linear Sampling method and the Range Test method. For notational simplicity we omit the dependance of the far-field operator and the solution of (\ref{esm}) on $\Omega$.
\begin{lemma}
Let the test domain $\Omega\supset D$ be a convex sound-soft obstacle with the $C^2$-smooth boundary $\partial \Omega$ and assume that $k^2$ is not the Dirichlet eigenvalue of $-\Delta$ over $\Omega$. If $\varphi \in L^2(\partial \Omega)$ is the solution of (\ref{range}) and $\tilde g \in L^2(\mathbb{S})$ is the solution of (\ref{esm}) with $B_{z,h}=\Omega$. Then we have
	\begin{equation}
		c_1\,|| \varphi ||_{L^2(\partial \Omega)} \leqslant||\tilde g||_{L^2(\mathbb{S})} \leqslant c_2\,|| \varphi ||_{L^2(\partial \Omega)},
	\end{equation}
	where $c_1, c_2$ are positive constants.
\end{lemma}
{\bf Proof.} By \cite[(1.55)]{Kirsch.2008}, the far-field operator $F=F_{\Omega}$ corresponding to the sound-soft obstacle $\Omega$ can be decomposed into the form
\begin{equation}\label{factorization-fhsh}
	F= -H^* S^{-1} H,\quad H: L^2(\s)\rightarrow H^{1/2}(\partial\Omega).
\end{equation}
Inserting \eqref{factorization-fhsh} into \eqref{esm} with $B_{z,h}$ replaced by $\Omega$ and combining \eqref{range} and \eqref{esm}, we get
\ben
F\tilde{g}=-H^* S^{-1} H\tilde g=H^*\varphi,\quad \varphi\in H^{-1/2}(\partial\Omega).
\enn
Since $k^2$ is not the Dirichlet eigenvalue of $-\Delta$ over $\Omega$, the operators $H^*$ and $S$ are both injective. Thus,
\begin{equation}
	S \varphi(x)= -H\tilde g(x), \quad x \in \partial \Omega.
\end{equation}
Note that the above equality is understood in $H^1(\partial \Omega)$, since
 $\partial \Omega$ is of $C^2$-smooth and $S: L^2(\partial \Omega)\rightarrow H^1(\partial \Omega)$ is bounded.
Then, for $\psi\in L^2(\Omega)$ and $\tilde{g}\in L^2(\s)$ we get
\begin{equation}
		\langle H^* \psi, \tilde{g} \rangle_{L^2(\mathbb{S})} = \langle \psi, H \tilde{g} \rangle_{L^2(\partial \Omega)}
		= -\langle \psi, S \varphi \rangle_{L^2(\partial \Omega)}
		= -\langle S^* \psi, \varphi \rangle_{L^2(\partial \Omega)}.
\end{equation}
Hence,
\begin{equation}\label{eq:6}
	\begin{split}
		||\tilde{g}||_{L^2(\mathbb{S})} &= \sup_{||H^* \psi||_{L^2(\s)} =1} \left|\langle H^* \psi, \tilde{g} \rangle_{L^2(\mathbb{S})}\right| \\
		&= \sup_{||H^* \psi||_{L^2(\s)} =1} \left| \langle S^* \psi, \varphi \rangle_{L^2(\partial \Omega)}\right| \\
		&\leqslant \sup_{||H^* \psi||_{L^2(\s)} =1} ||S^* \psi||_{L^2(\partial \Omega)} \  ||\varphi||_{L^2(\partial \Omega)} \\
		&\leqslant \sup_{||H^* \psi||_{L^2(\s)} =1} c_4||\psi||_{H^{-1}(\partial \Omega)} \  ||\varphi||_{L^2(\partial \Omega)}\\
&\leqslant \sup_{||H^* \psi||_{L^2(\s)} =1} c_4||\psi||_{L^{2}(\partial \Omega)} \  ||\varphi||_{L^2(\partial \Omega)}\\
		&\leqslant c_2  ||\varphi||_{L^2(\partial \Omega)}\\
\end{split}
\end{equation}

On the other hand, using again the assumption on $k^2$ and the smoothness of $\partial \Omega$, for every $d\in\s$ we can always find a $\psi(\cdot;d)\in L^2(\partial \Omega)$ such that the equality
\ben
-e^{ikx\cdot d}=[S\psi(\cdot;d)](x),\quad x\in \partial \Omega,
\enn
holds in the sense of $H^1(\partial \Omega)$.
Moreover, it holds that $||\psi(\cdot;d)||_{L^2(\partial \Omega)}\leq C$ uniformly in all $d\in \s$.
Hence,
\begin{equation}\label{far-G}
u_{\Omega}^{\infty}(\hat{x};d)=\int_{\partial \Omega} e^{-i k \hat{x} \cdot y} \psi(y;d) d s(y), \quad\hat{x} \in \mathbb{S}.
\end{equation}
Since $H^* \varphi = u^{\infty}_{D}= F\tilde{g}$, using \eqref{far-G} and \eqref{range-int} we get
\begin{equation}
	\begin{split}
		\int_{\partial \Omega} e^{-i k \hat{x} \cdot y} \varphi(y) d s(y) &=\int_{\mathbb{S}} \int_{\partial \Omega} e^{-i k \hat{x} \cdot y} \psi(y ; d) d s(y) \tilde{g}(d) d s(d) \\
		&=\int_{\partial \Omega} e^{-i k \hat{x} \cdot y}\left(\int_{\mathbb{S}} \psi(y ; d) \tilde{g}(d) d s(d)\right) d s(y),
	\end{split}
\end{equation}
implying that
\begin{equation}
	\int_{\partial \Omega} e^{-i k \hat{x} \cdot y}\left(\varphi(y) - \int_{\mathbb{S}} \psi(y ; d) \tilde{g}(d) d s(d)\right) d s(y) = 0, \quad\forall\;\hat{x}\in \s.
\end{equation}
By the injectivity of $H^*$,
\begin{equation}
	\varphi(y) = \int_{\mathbb{S}} \psi(y ; d) \tilde{g}(d) d s(d).
\end{equation}
Thus,
\begin{equation}
	\begin{split}
		|| \varphi ||^2_{L^2(\partial \Omega)} &= \int_{\partial \Omega} \left|  \int_{\mathbb{S}} \psi(y ; d) \tilde{g}(d) d s(d) \right|^2 ds(y) \\
		&\leqslant \int_{\partial \Omega} ||\psi(y;\cdot)||^2_{L^2(\mathbb{S})} ds(y) ||\tilde{g}||^2_{L^2(\mathbb{S})} \\
		&\leqslant \frac{1}{c_1} ||\tilde{g}||^2_{L^2(\mathbb{S})}.
	\end{split}
\end{equation}
This together with \eqref{eq:6} proves the relaiton $c_1|| \varphi ||_{L^2(\partial \Omega)} \leqslant||\tilde{g}||_{L^2(\mathbb{S})} \leqslant c_2|| \varphi ||_{L^2(\partial \Omega)} $.  $\hfill\square$

\section{Numerical tests}\label{sec:4}
\subsection{Reconstruction of a finite number of point-like obstacles}

\begin{figure}[H]
	\centering    
	\subfigure 
	{
		\begin{minipage}[t]{0.45\linewidth}
			\centering          
		\includegraphics[scale=0.2]{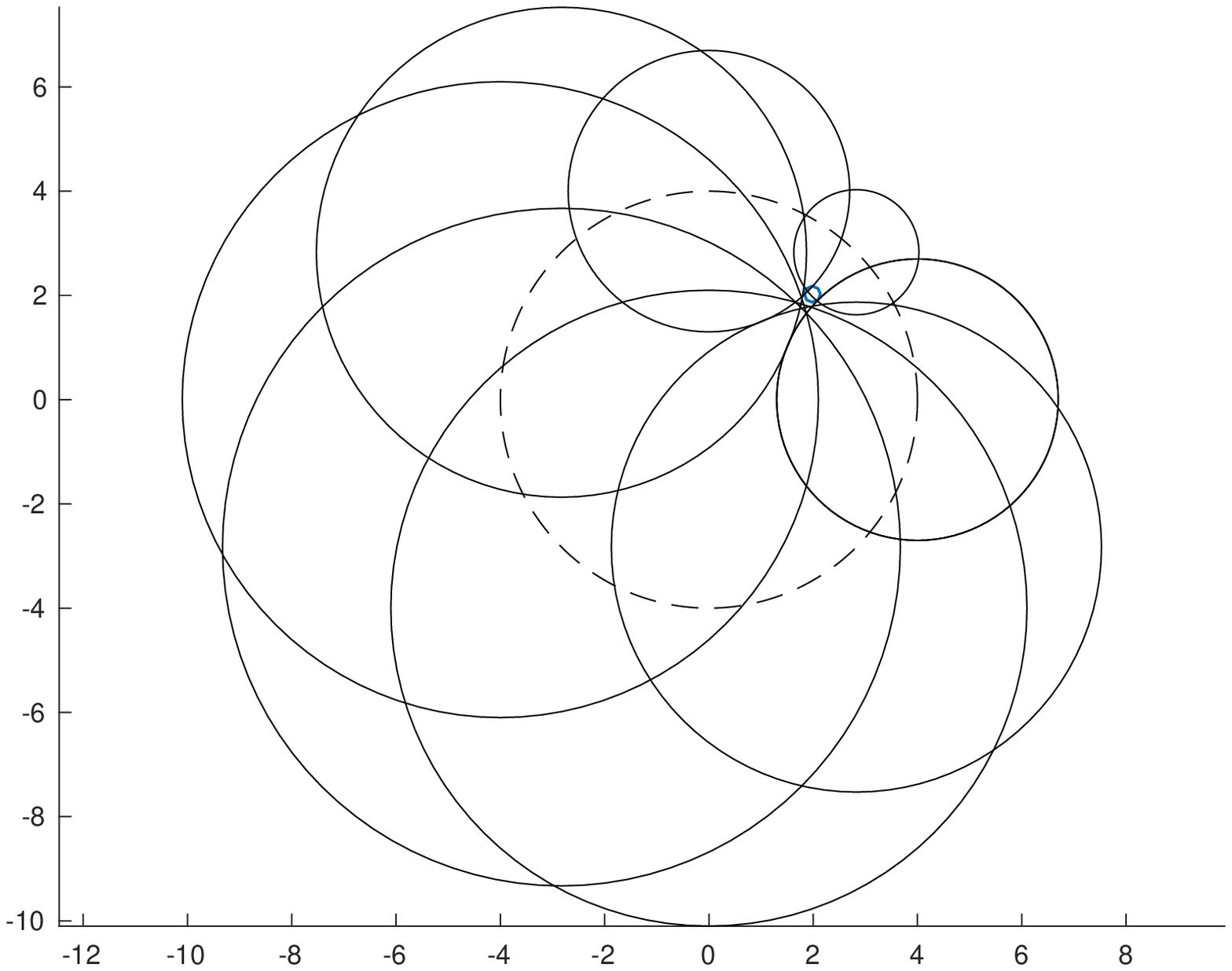}
			\caption*{$z^{*}=[2,2]$}
		\end{minipage}
	}	
	\subfigure 
	{
		\begin{minipage}[t]{0.45\linewidth}
			\centering      
			\includegraphics[scale=0.2]{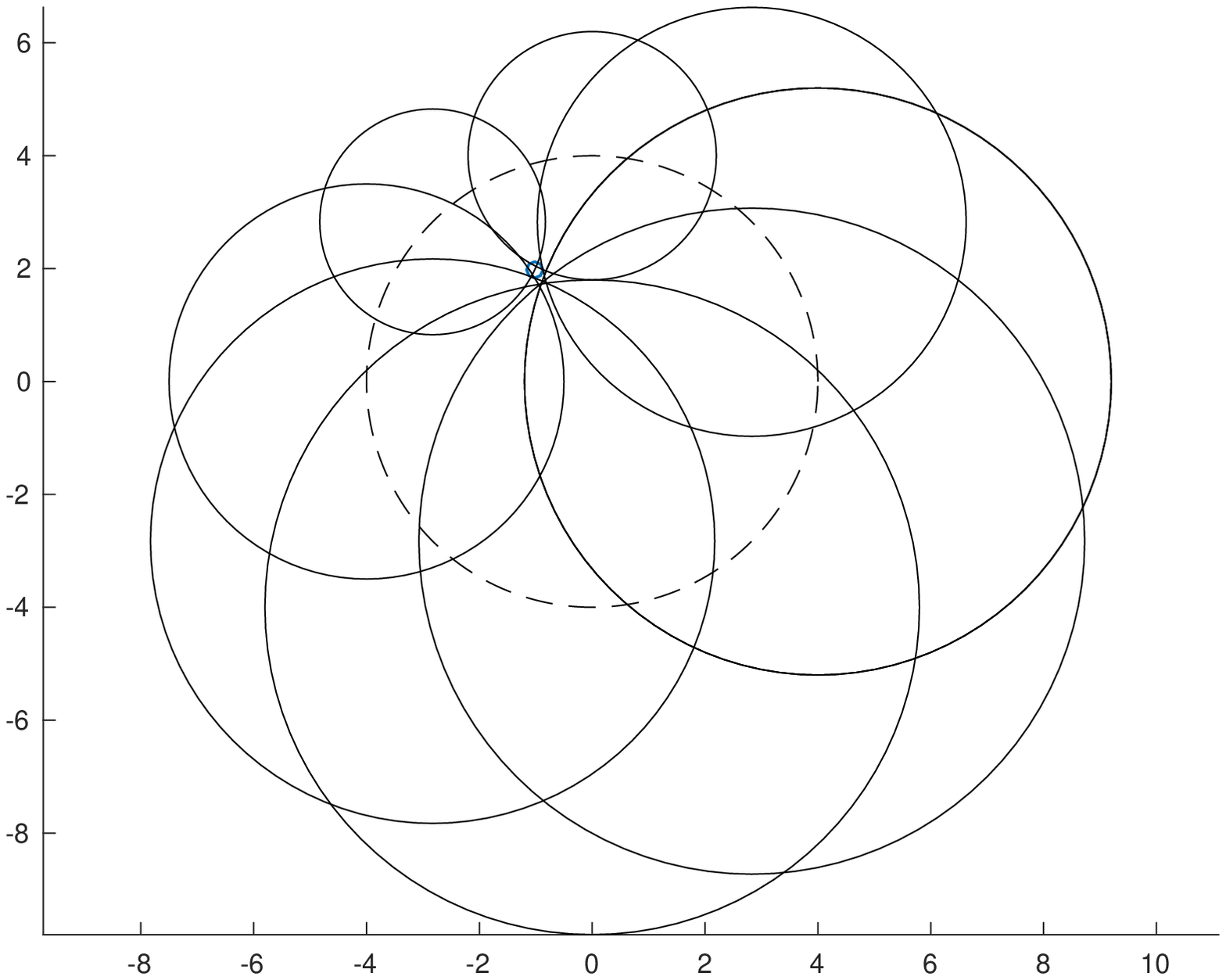}
			\caption*{$z^{*}=[-1,2]$}
		\end{minipage}
	}
	\subfigure 
	{
		\begin{minipage}[t]{0.45\linewidth}
			\centering      
			\includegraphics[scale=0.2]{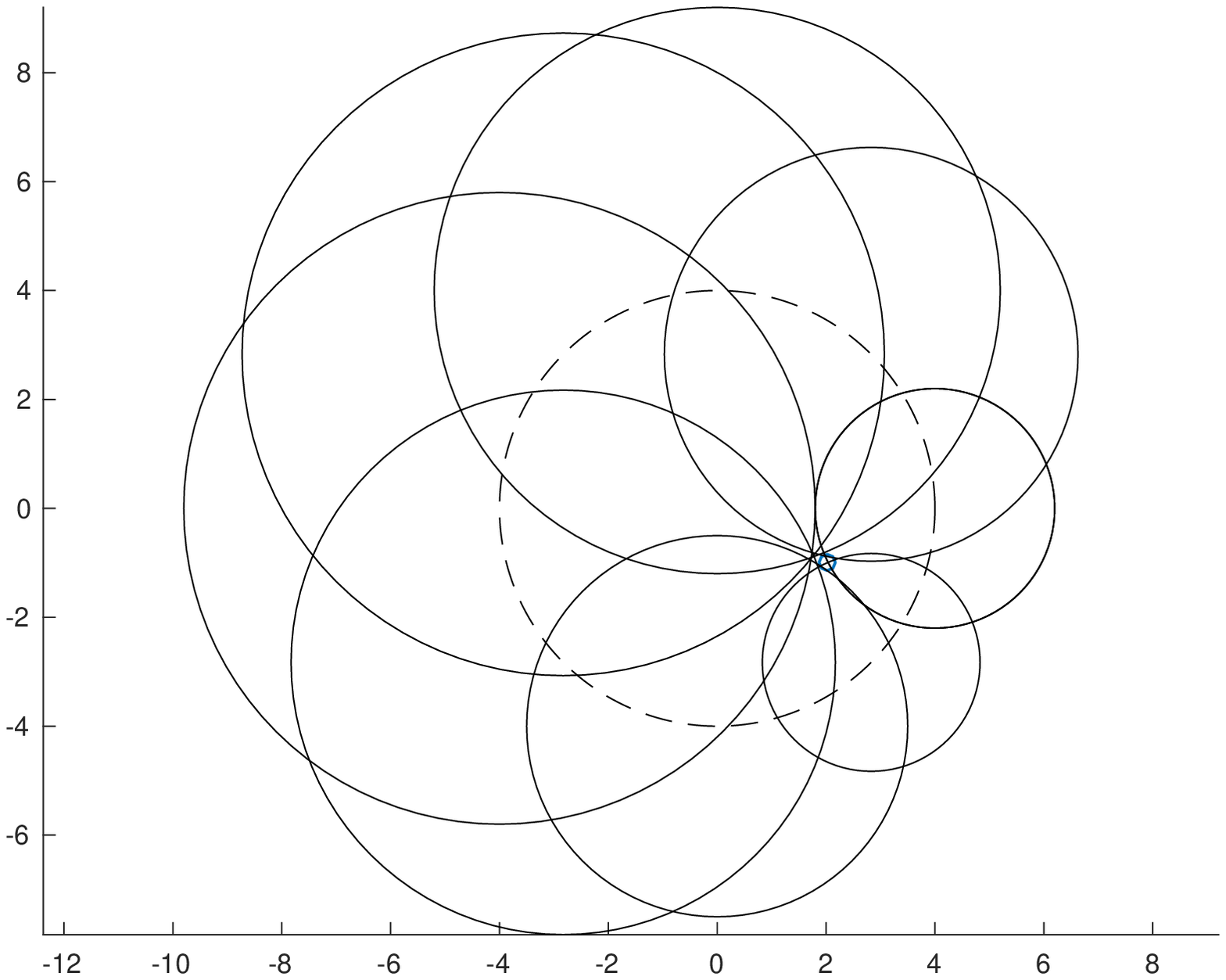}
			\caption*{$z^{*}=[2,-1]$}
		\end{minipage}
	}
	\subfigure
{
	\begin{minipage}[t]{0.45\linewidth}
		\centering
		\includegraphics[scale=0.2]{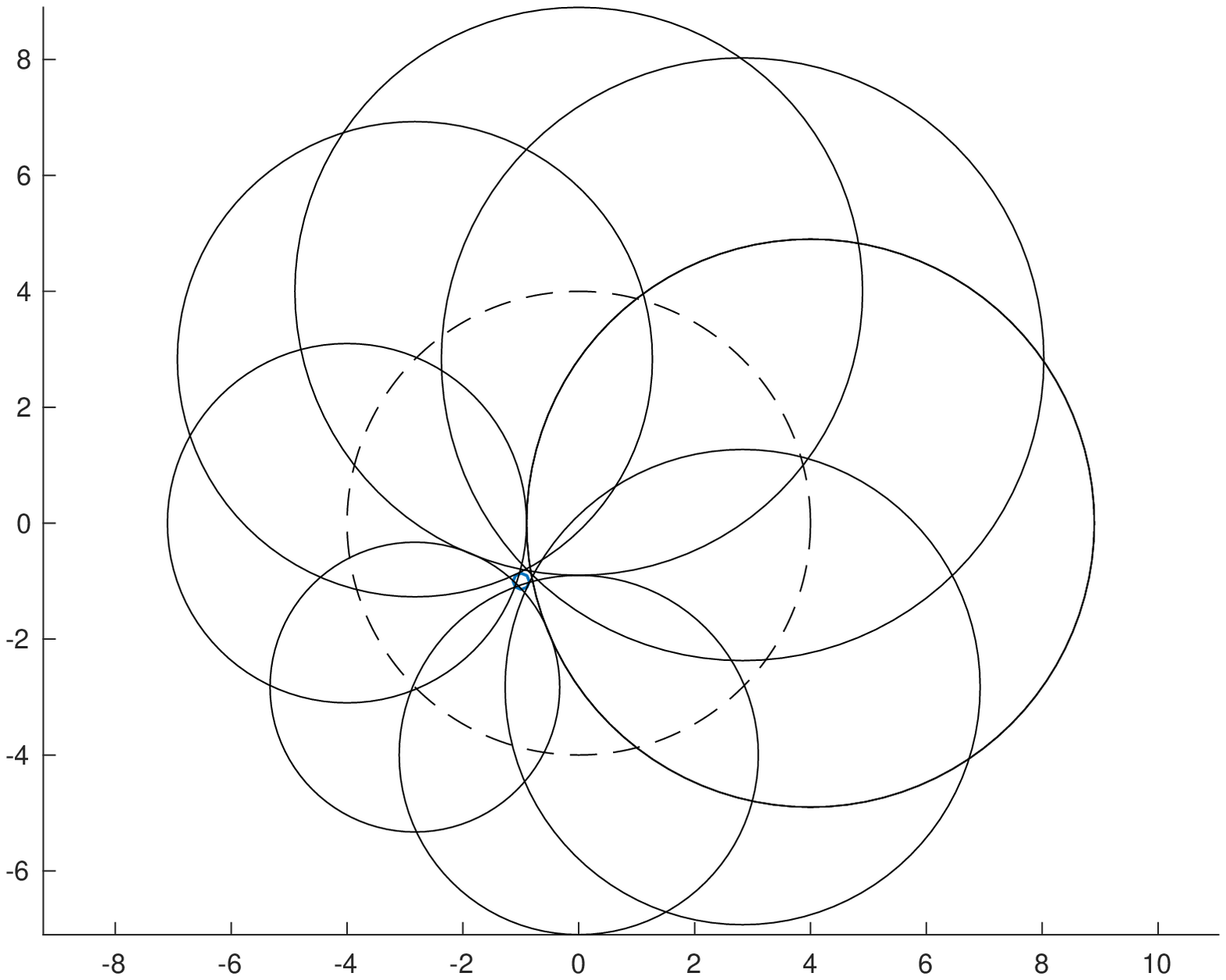}
		\caption*{$z^{*}=[-1,-1]$}
	\end{minipage}
}
	\caption{Imaging scheme I: locating one point-like obstacle/source from a single far-field pattern.}
	\label{point-like}
\end{figure}

\begin{figure}[H]
	\centering
	\includegraphics[width=0.4\linewidth]{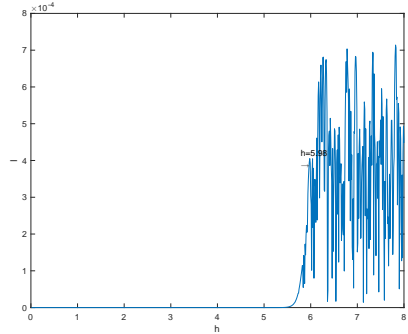}
	\caption{Figure of the function $h\rightarrow \widetilde{W}(z,h)$ with $z=[4,0]$, $h\in(0,8)$ for locating $z^*=[-2,0]$.}
	\label{H-rate}
\end{figure}

\begin{figure}[H]
	\centering    
	
	\subfigure 
	{
		\begin{minipage}[t]{0.3\linewidth}
			\centering          
			\includegraphics[scale=0.3]{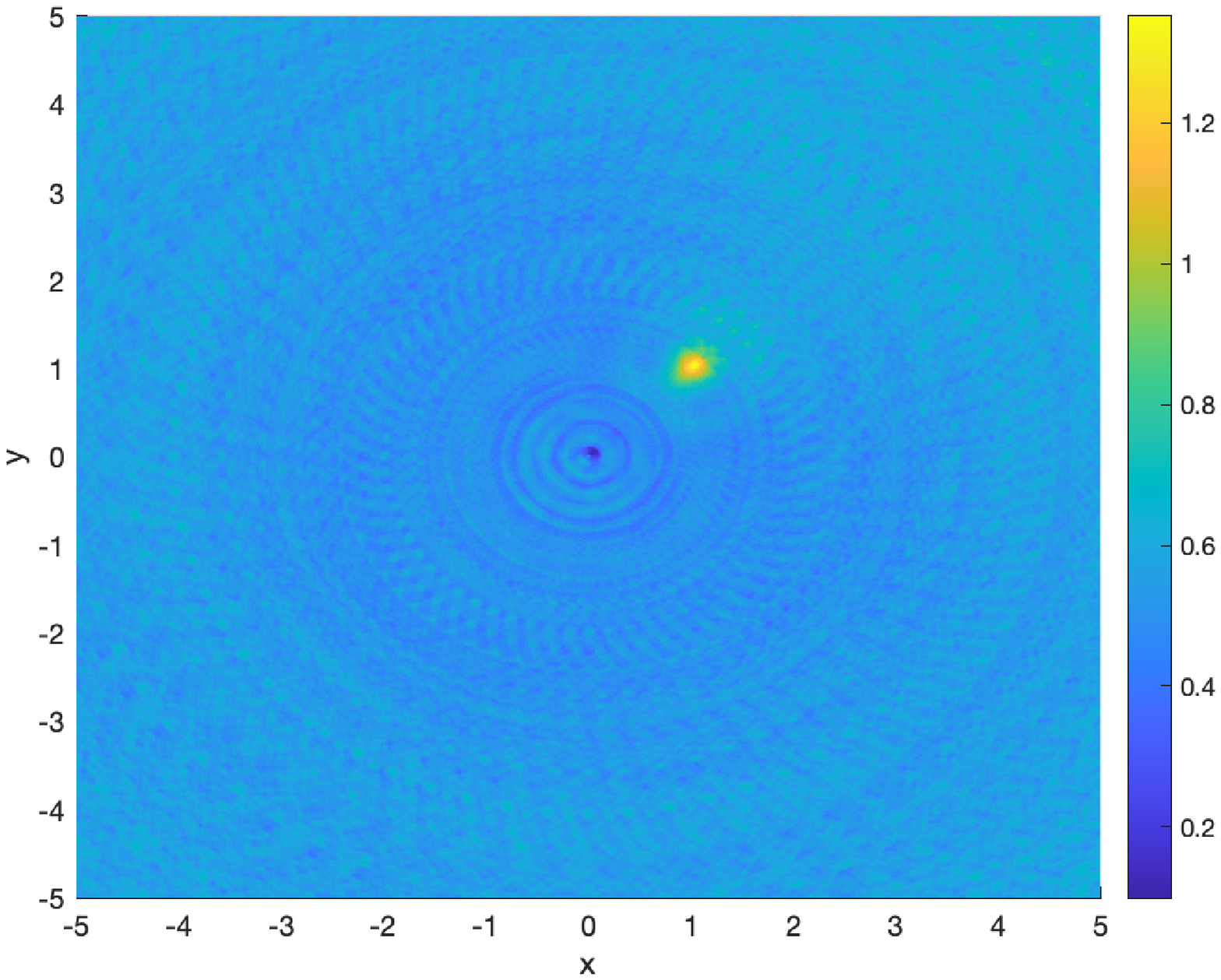}
		\end{minipage}
	}	
	\subfigure 
	{
		\begin{minipage}[t]{0.3\linewidth}
			\centering      
			\includegraphics[scale=0.3]{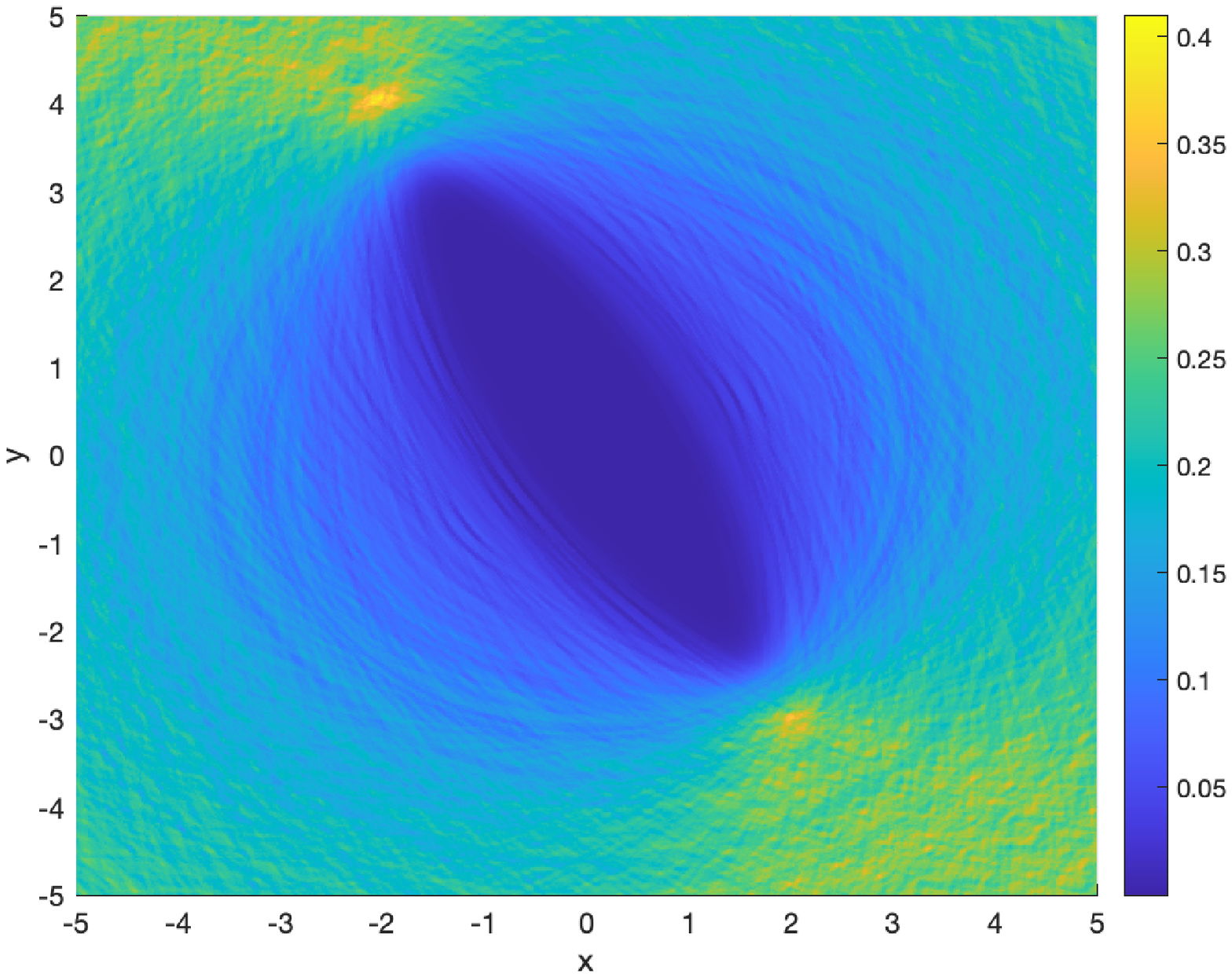}
		\end{minipage}
	}
	\subfigure 
	{
		\begin{minipage}[t]{0.3\linewidth}
			\centering      
			\includegraphics[scale=0.3]{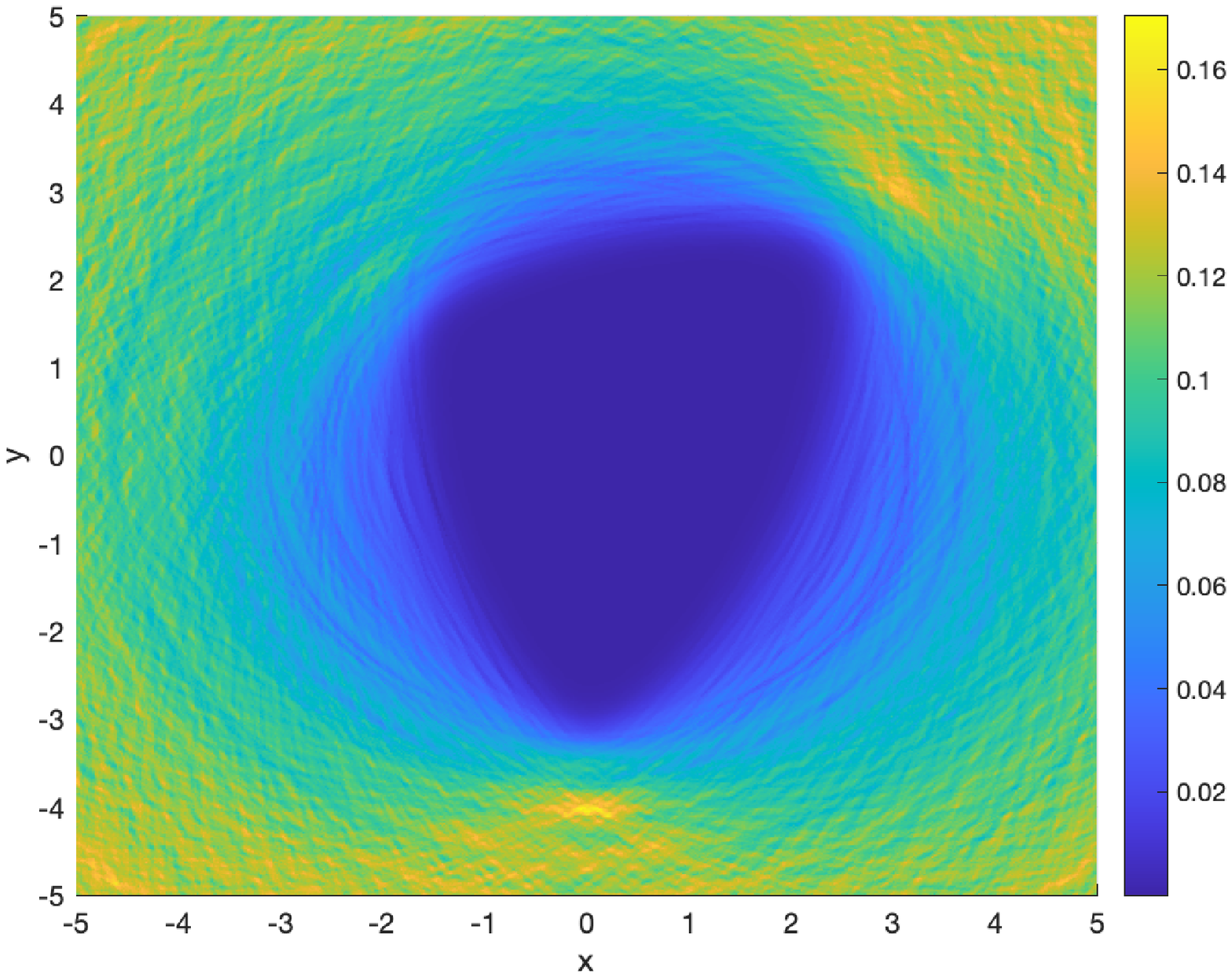}
		\end{minipage}
	}
	\caption{Reconstruction of one (left), two (middle) and three (right) points by applying imaging scheme II. The dotted circle denotes $\Gamma_R$.}
	\label{muti_point-like}
\end{figure}

%
Assume that the obstacle (resp. source support) $D$ has been shrank to a point located at $z^*$ (for instance, if the wavelength $\lambda=2\pi/k$ is much bigger than the diameter of $D$).
In this case, the scattered (resp. radiated) wave field can be asymptotically written as $c^*\Phi(x,z^*)$, where $c^*\in \C$ depends on the incoming wave, the scattering strength of $D$ as well as the location point $z^*$. For simplicity we suppose that $c^*=1$. Then the far-field pattern
 $u^{\infty}=u_D^{\infty}$ takes the simple form
 \begin{equation}
 	u^{\infty}(\hat x) = e^{-ikz^{*} \cdot \hat{x}}.
 \end{equation}
To perform numerical examples, we set the wave number $k=6$ and suppose that $z^*\in B_R$ with $R=4$. The number of sampling centers $z_n$ lying on $|x|=4$ is taken to be $N_z=8$, and the parameter for truncating the infinite series (\ref{Wzh}) is chosen as $N=60$.
 The threshhold specified in our imaging scheme is set as $\delta=4\times 10^{-4}$. In these settings
we obtain Figure \ref{point-like} where $z^*=[2,2], [-1,2], [2,-1], [-1,-1]$ can be accurately located with a single far-field pattern without polluted noise.
In Fig. \ref{point-like}, the dotted curve represents the circle $\Gamma_{R}$ where the centers of the test disks are located. The solid circles are boundaries of the test disks $B_{h_{z_n}}(z_n)$, where $h_{z_{n}}$ denotes the distance between $z_n$ and $z^*$. In Figure \ref{H-rate}, we plot the function $h\rightarrow \widetilde{W}(z,h)$ with $z=[4,0]$ for locating $z^*=[-2,0]$. Obviously, we have $|z-z^*|=6$. It is seen from Figure \ref{H-rate} that the values of the function $h\rightarrow \widetilde{W}(z,h)$ for $h\in(0,6)$ are much smaller than those for $h\in(6,8)$.
With our threshhold the distance between $z$ and $z^*$ is calculated as 5.98.
Although the first scheme can be used previously to locate a point-like obstacle/source, it is not straightforward for imaging an extended obstacle/source.
Below we describe a more direct imaging scheme II.
\begin{itemize}
	\item Suppose that $B_R \supset D$ for some $R>0$ and collect the measurement data $u_D^{\infty}(\hat{x})$ for all $\hat{x} \in \mathbb{S}$. Let $Q\supset D$ be our search/computational region for imaging $D$;
	\item Choose sampling centers $z_{n} \in \Gamma_{R}:= \{x:\ |x|=R\}$ for $n = 1,...,N_{z}$ and
	choose sampling radii $h_m\in(0,2R)$ to get different spectral systems $(\lambda^{(j)}_{z_n,h_m}, \varphi^{(j)}_{z_n,h_m})$ (see (\ref{eigsys-sf}) or (\ref{eigsys}));
	\item For each $z_n\in \Gamma_R$, define the function $\mathcal{I}_{n}(x)=\widetilde{W}(z_n,|x-z_n|)$ for $x\in Q$ (see \eqref{Wt});
\item The imaging function for recovering $D$ is defined as $\mathcal{I}(x)= \sum_{n=1}^{N_z}\mathcal{I}_{n}(x)$, $x\in Q$;
\end{itemize}
Note that, by Corollary \ref{W}, $\mathcal{I}_n(x)>0$ if $|x-z_n|>\max_{y\in D} |z_n-y|$ and $\mathcal{I}_n(x)=0$ if otherwise. Hence, the values of $\mathcal{I}(x)$ for $x\in Q\backslash\overline{D}$ should be larger than those for $x\in D$.
As an example, we apply this new scheme to image
$D=\bigcup_{j=1,2,\cdots,M}\{z_j^*\}$ which consists of multiple point-like scatterers. For simplicity we neglect the multiple scattering between them and write the far-field pattern as
 \begin{equation}
u^{\infty}(z^{*}) = e^{-ikz^{*}_1 \cdot \hat{x}}+e^{-ikz^{*}_2 \cdot \hat{x}}+...+e^{-ikz^{*}_M \cdot \hat{x}}.
\end{equation}
The results for
reconstructing one point $z_1^{*}=[1,1]$ ($M=1$),  two points $z_1^{*}=[-2,4]$, $z_2^*=[2,-3]$ ($M=2$) and three points $z_1^{*}=[3,3]$, $z_2^*=[-2,2]$, $z_3^*=[0,-4]$ ($M=3$) are shown in
Fig \ref{muti_point-like}. As one can imagine, with a single far-field pattern our approach can only  recover the convex hull of these points. In the case of two points, the line segment connecting $z_1^*$ and $z_2^*$ is shown in the middle of Fig. \ref{muti_point-like}. The triangle formed by $z_1^*$, $z_2^*$ and $z_3^*$ is shown in the right figure.

\begin{remark}
It is important to remark that there exist other approaches for locating a finite number of point scatterers using several incident waves, for example the MUSIC algorithm (which can be regarded as the discrete analogue of the classical factorization method \cite{kirsch2002}). Our concern here is to show the capability of the one-wave factorization method for identifying singularities of the scattered/radiated wave field. In 2D, the singularity is of 
logarithmic type at the point-like scatterers.
\end{remark}

%
\subsection{Reconstruction of a triangular source support}
Suppose that the source support $D$ is a triangle with the three corners located at $(-2,-2)$,$(-2,2)$ and (2,-2). The source function is supposed to be a constant. For simplicity we assume that $\chi_D f(x)\equiv 1$ on $\overline{D}$, so that the far-field patten takes the explicit form
\begin{equation}
u^{\infty}(\hat{x}) = \frac{i }{4} \int_{D} e^{-i k \hat{x} \cdot z} dz
                           = \frac{i}{4} \int_{-2}^{2} e^{-i k \hat{x}_1z_1} \int_{-2}^{-z_1} e^{-i k \hat{x}_2z_2} dz_2 dz_1,
\end{equation}
where $\hat{x}=(\hat{x}_1,\hat{x}_2)\in \s$ and $z=(z_1, z_2)\in \R^2$.
We want to test the sensitivity of our approach to the incident wavenumber $k$ and to the circle $\Gamma_R=\{x: |x|=R\}=\partial B_R$. In this subsection, the number of sampling centers $z_n$ lying on $|x|=R$ is set to be $N_z=64$ and the truncation parameter to be $N=80$. 

\begin{figure}[H]
	\centering
	\subfigure
	{
		\begin{minipage}[t]{0.3\linewidth}
			\centering
			\includegraphics[scale=0.3]{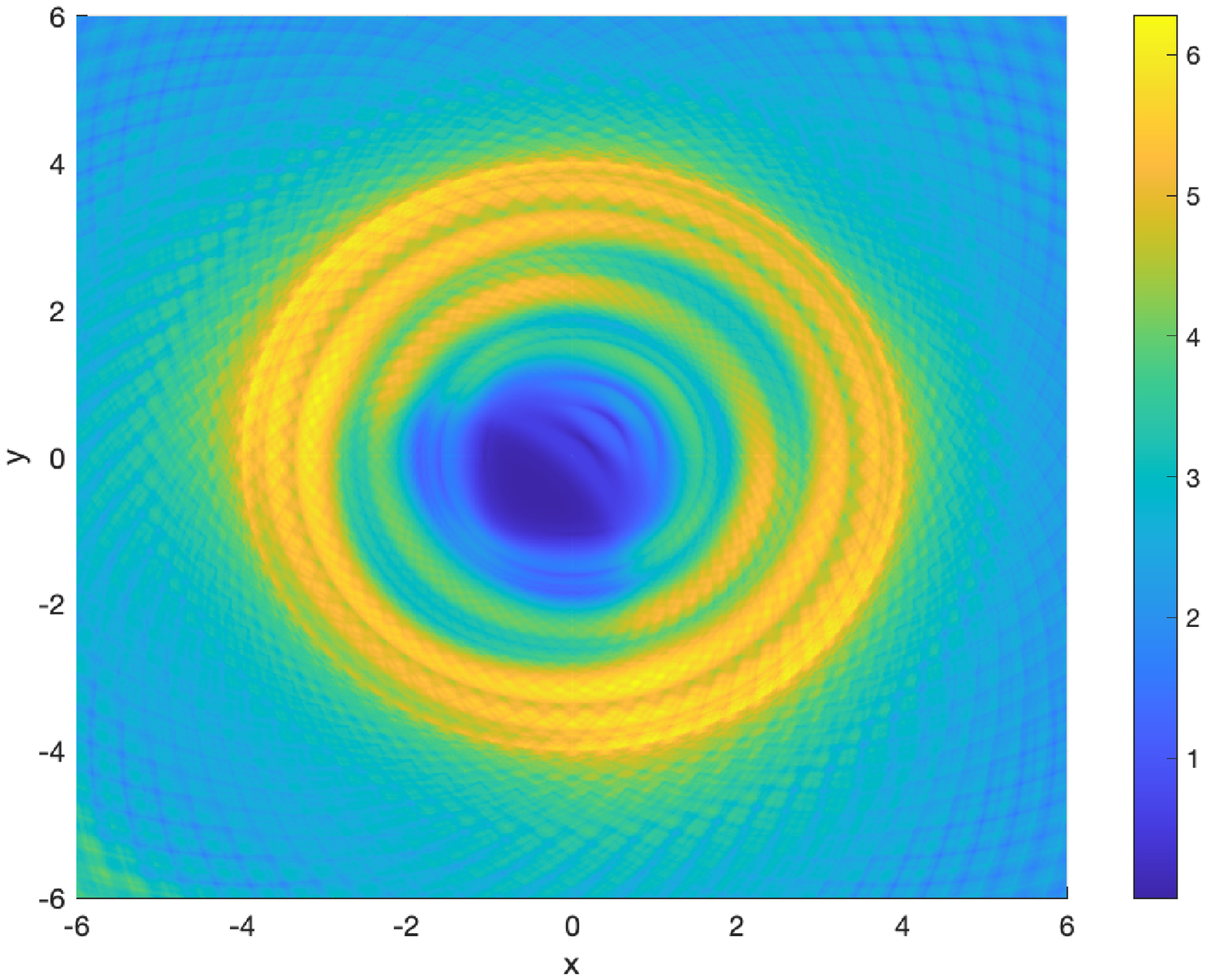}
			\caption*{$k=1.5$}
		\end{minipage}
	}
	\subfigure
	{
		\begin{minipage}[t]{0.3\linewidth}
			\centering
			\includegraphics[scale=0.3]{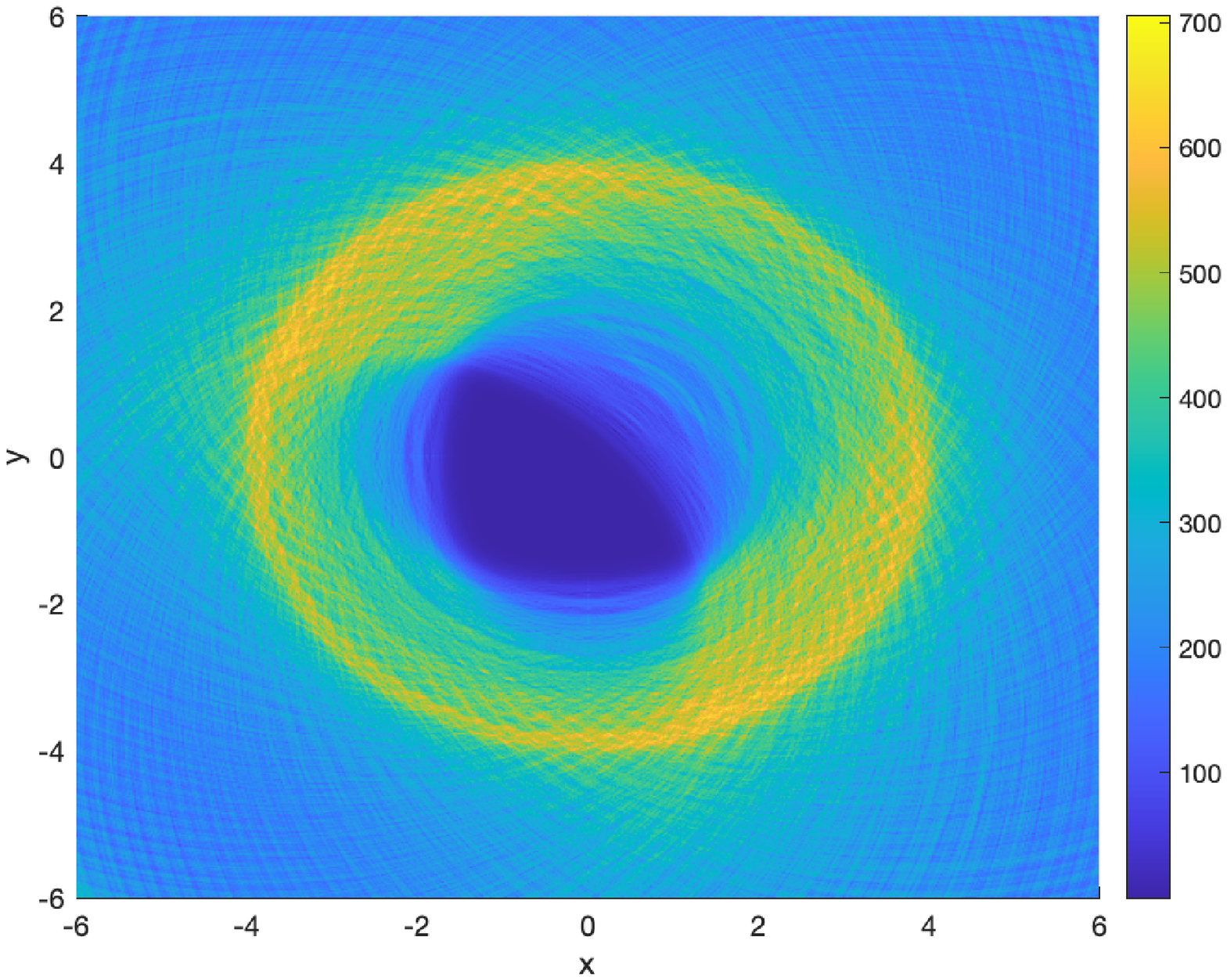}
			\caption*{$k=6$}
		\end{minipage}
	}
	\subfigure
{
	\begin{minipage}[t]{0.3\linewidth}
		\centering
		\includegraphics[scale=0.3]{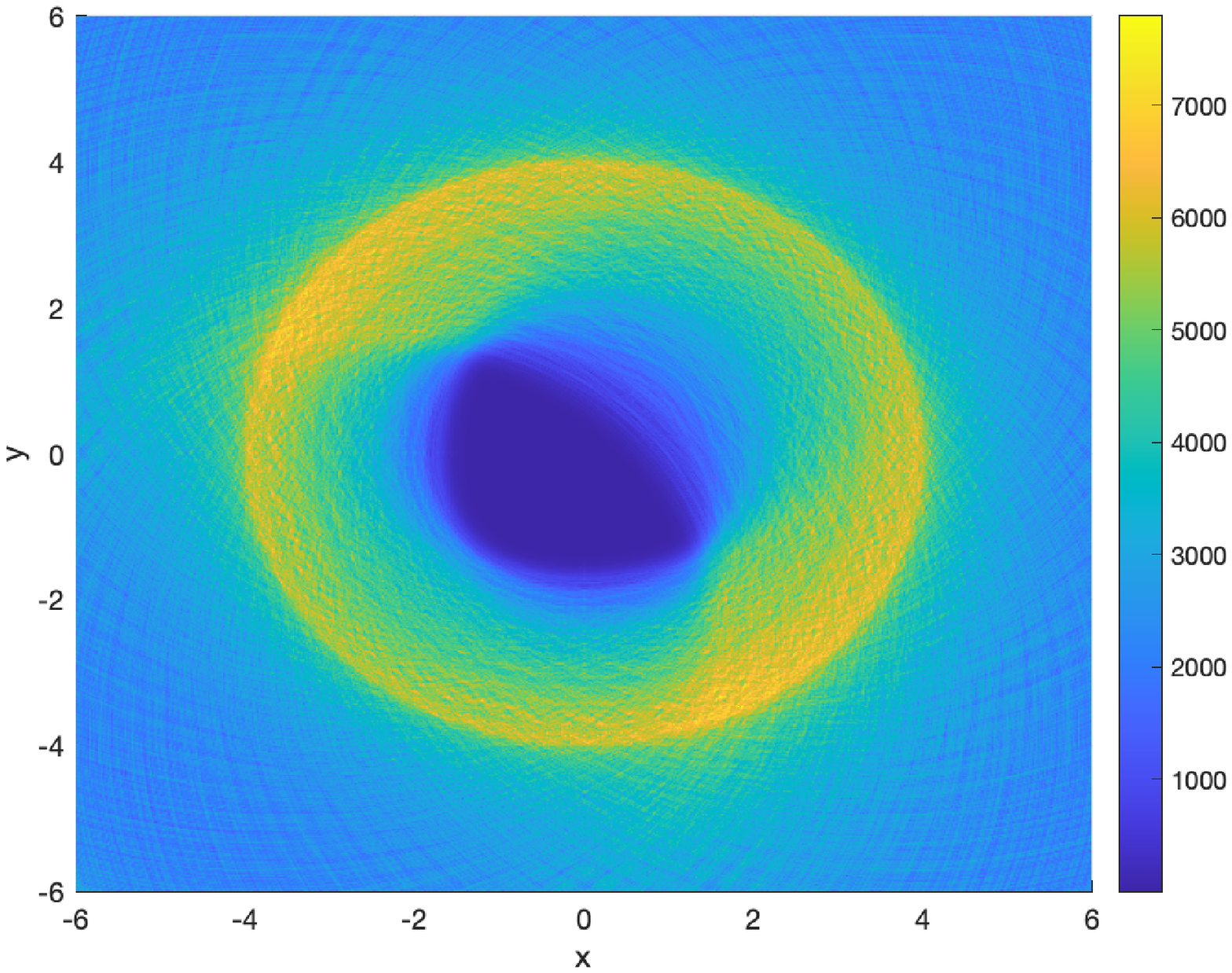}
		\caption*{$k=12$}
	\end{minipage}
}
	\caption{Image of a triangular source support from a single far-field pattern excited at different energies.}
	\label{Source with different k}
\end{figure}

Firstly, we fixed $R=4$ and change the excited frequencies. The far-field data corresponding to $B_h(z)$ are supposed to be excited at the same frequency as for $D$.
 In Fig. \ref{Source with different k}, the recovery of $D$ from a single far-field pattern at different frequencies are illustrated. We observe that a regularization parameter depending on $k$ must be properly selected to get a satisfactory image of $D$. In our numerical tests, $\alpha>0$ is chosen by the method of trial and error. Consequently,
 we take  $\alpha=1e-22$, $1e-13$ and $1e-8$ corresponding to the wavenumbers
 $k=1.5, 6, 12$, respectively. From Fig. \ref{Source with different k}, one can conclude that
 a better image can be achieved at higher frequencies.

\begin{figure}[H]
	\centering
	\subfigure
	{
		\begin{minipage}[t]{0.3\linewidth}
			\centering
			\includegraphics[scale=0.3]{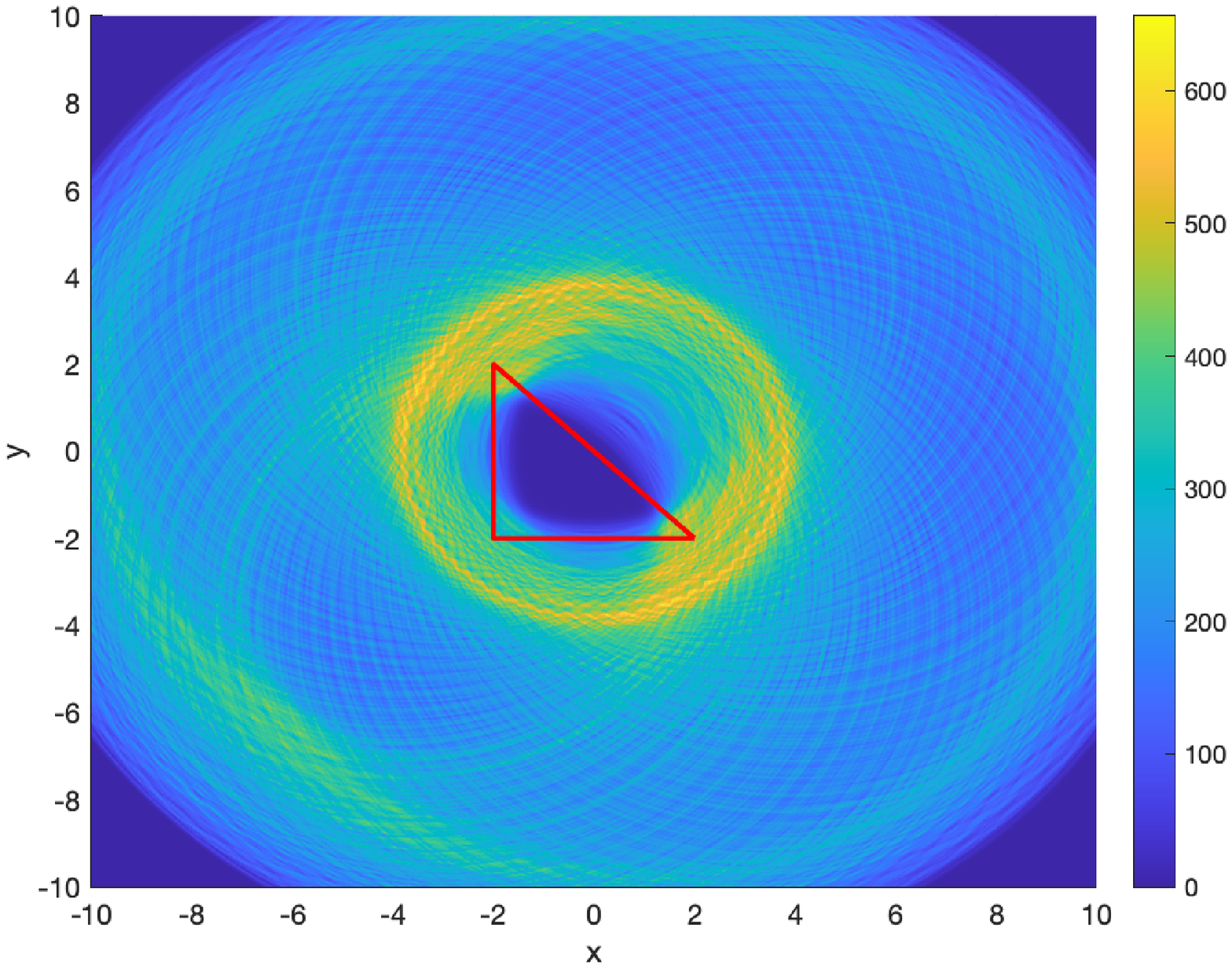}
			\caption*{$R=4$}
		\end{minipage}
	}
	\subfigure
	{
		\begin{minipage}[t]{0.3\linewidth}
			\centering
			\includegraphics[scale=0.3]{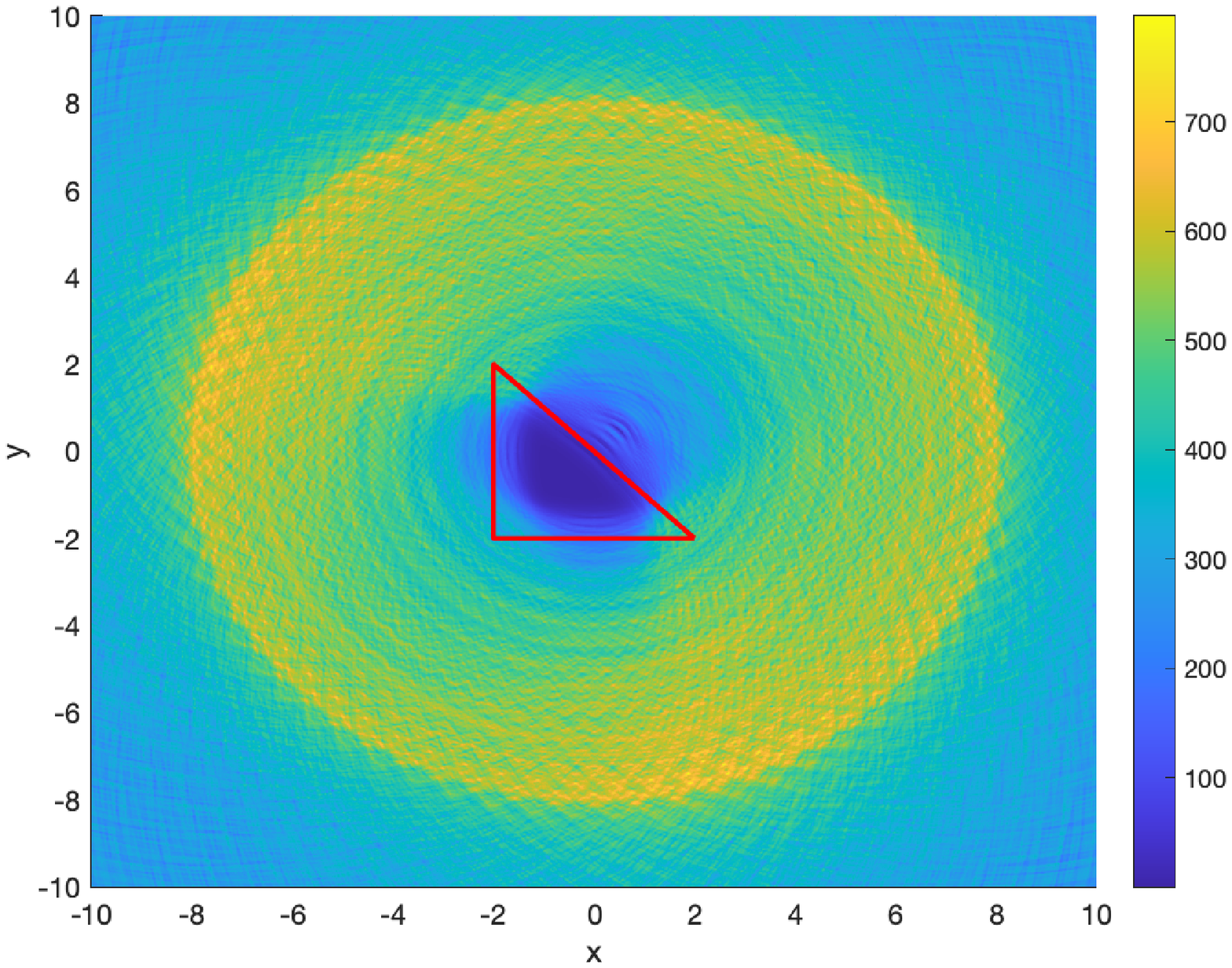}
			\caption*{$R=8$}
		\end{minipage}
	}
	\subfigure
{
	\begin{minipage}[t]{0.3\linewidth}
		\centering
		\includegraphics[scale=0.3]{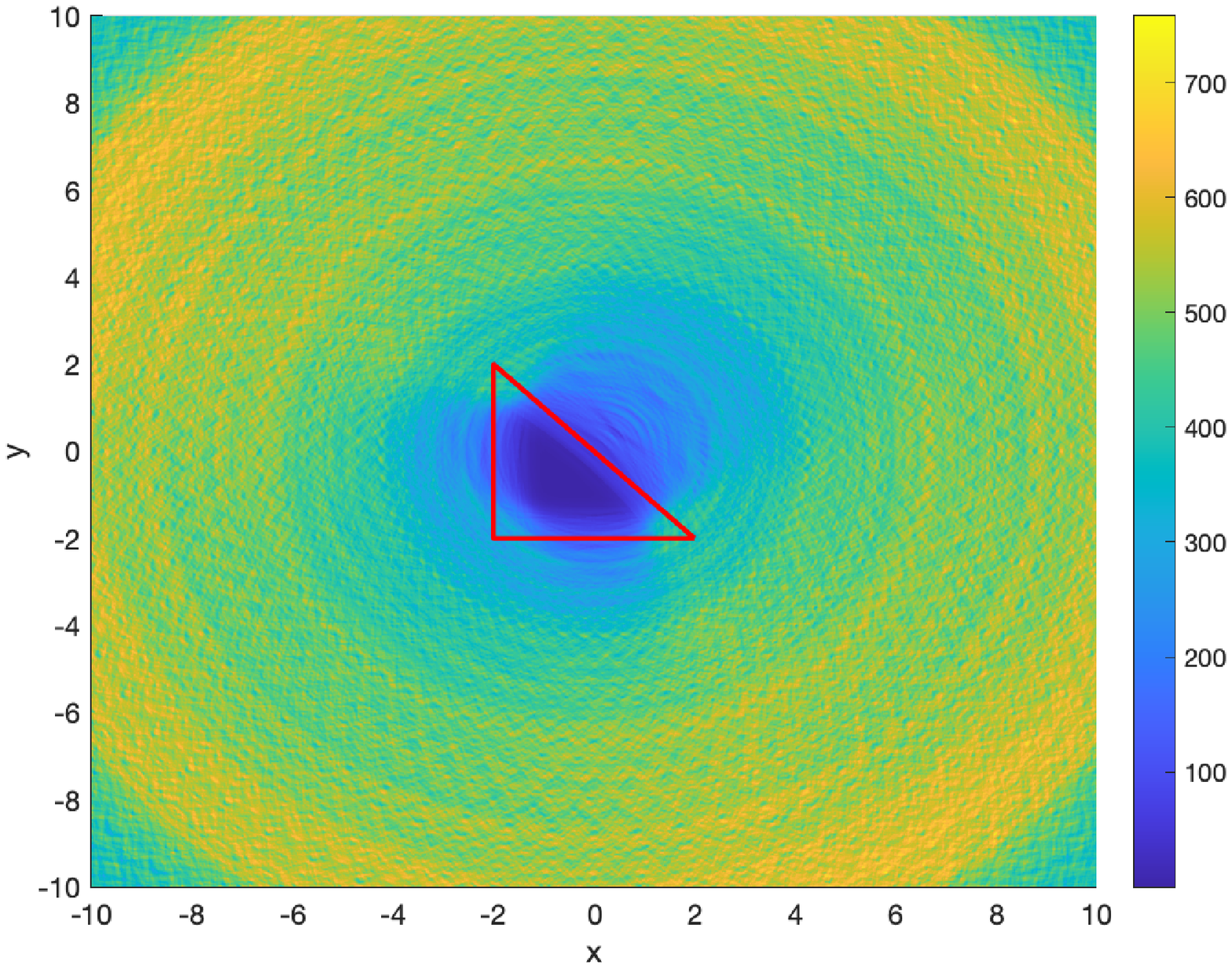}
		\caption*{$R=12$}
	\end{minipage}
}
	\caption{Image of a triangular source support $D$ by using the spectral system of $F_{z,h}$ for $z\in \Gamma_R$. We fix $k=6$.}
	\label{Source using different l}
\end{figure}

 Secondly, we fix $k=6$ and recovery $D$ by using sampling disks with the centers equally distributed on $\Gamma_R$ with $R=4, 8, 12$. The regularization parameter are chosen as $\alpha=1e-13$, $1e-13$ and $1e-20$. Note that a smaller $R$ gives more a priori information on $D$. It is seen from Figure \ref{Source using different l} that a larger $R$ yields a worse image of $D$.

\subsection{Reconstruction of a sound-soft square obstacle}
\begin{figure}[H]
		\centering
		\subfigure
		{
		\begin{minipage}[t]{0.45\linewidth}
			\centering \includegraphics[scale=0.4]{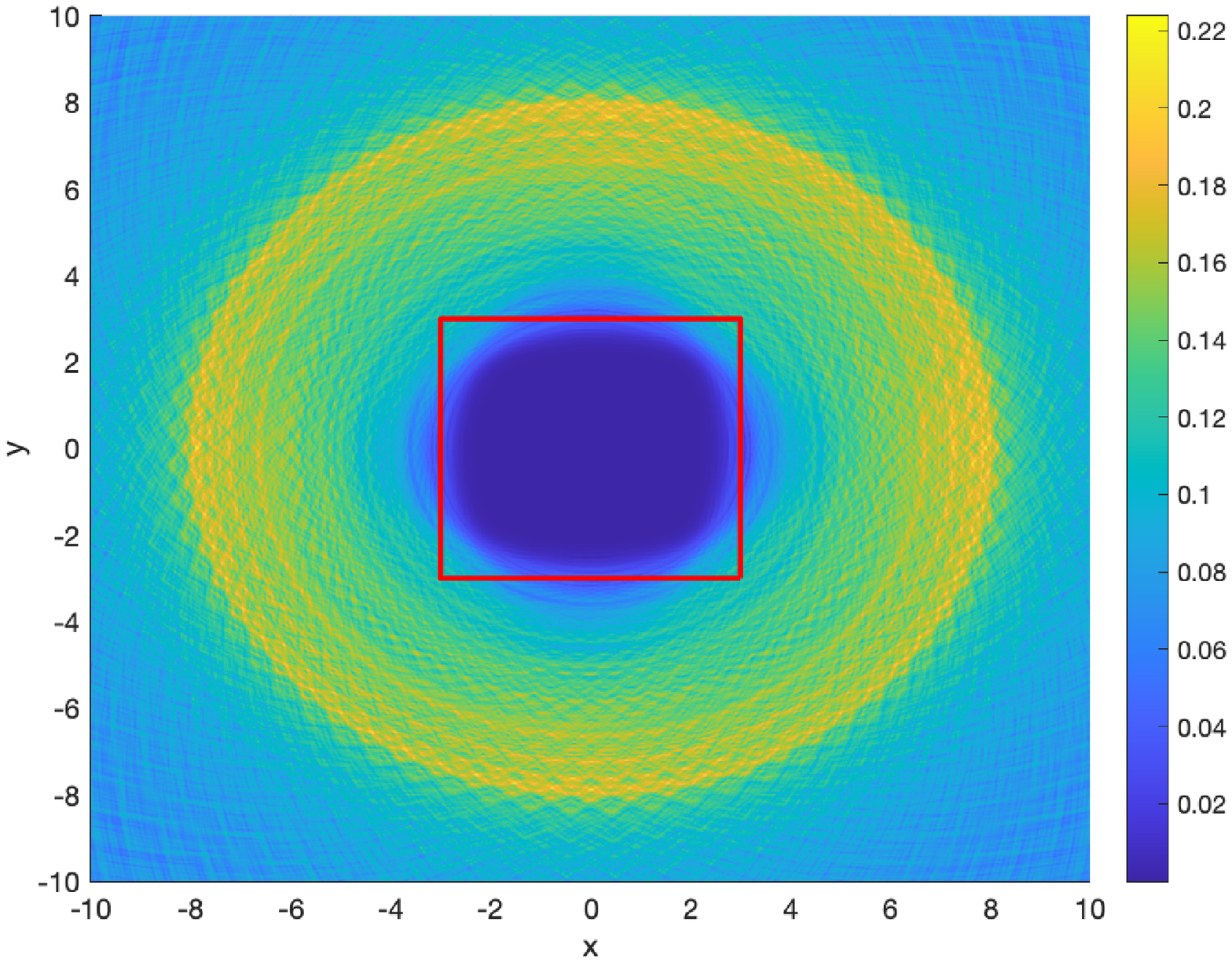}
    	\end{minipage}
        }
		\subfigure
{
	\begin{minipage}[t]{0.45\linewidth}
		\centering \includegraphics[scale=0.4]{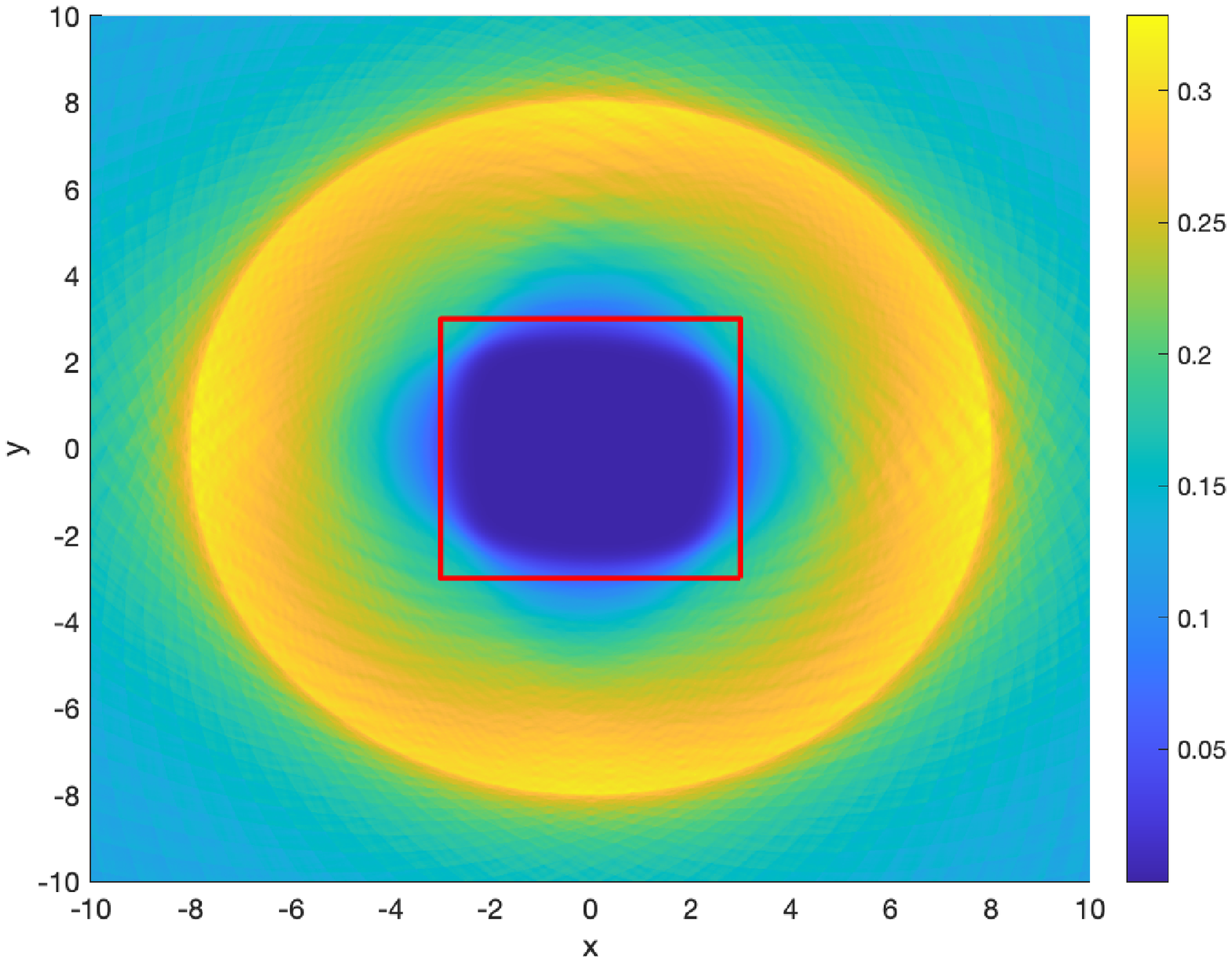}
	\end{minipage}
}
	\caption{Reconstruction of a sound-soft square $D=[-3,3]\times[-3,3]$ by using test/sampling disks with different boundary conditions.
Left: $B_h(z)$ are sound-soft; Right: $B_h(z)$ are of impedance type.}
    \label{Obstacle using different boundary condition test disks}
\end{figure}

\begin{figure}[H]
			\centering
    		\subfigure
{
	\begin{minipage}[t]{0.3\linewidth}
		\centering
		\includegraphics[scale=0.3]{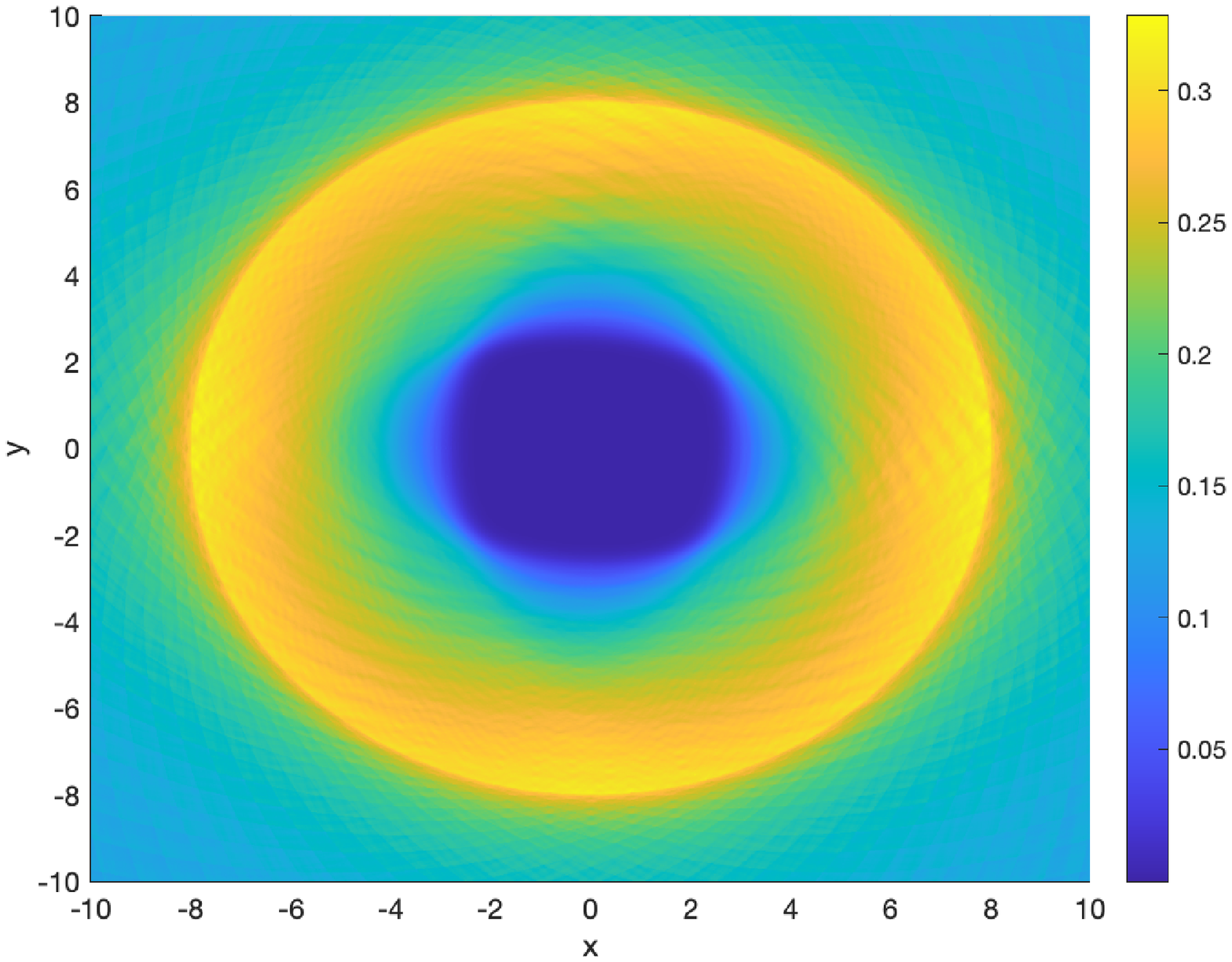}
		\caption*{$\theta=2$ rad}
	\end{minipage}
}
		\subfigure
{
	\begin{minipage}[t]{0.3\linewidth}
		\centering
		\includegraphics[scale=0.3]{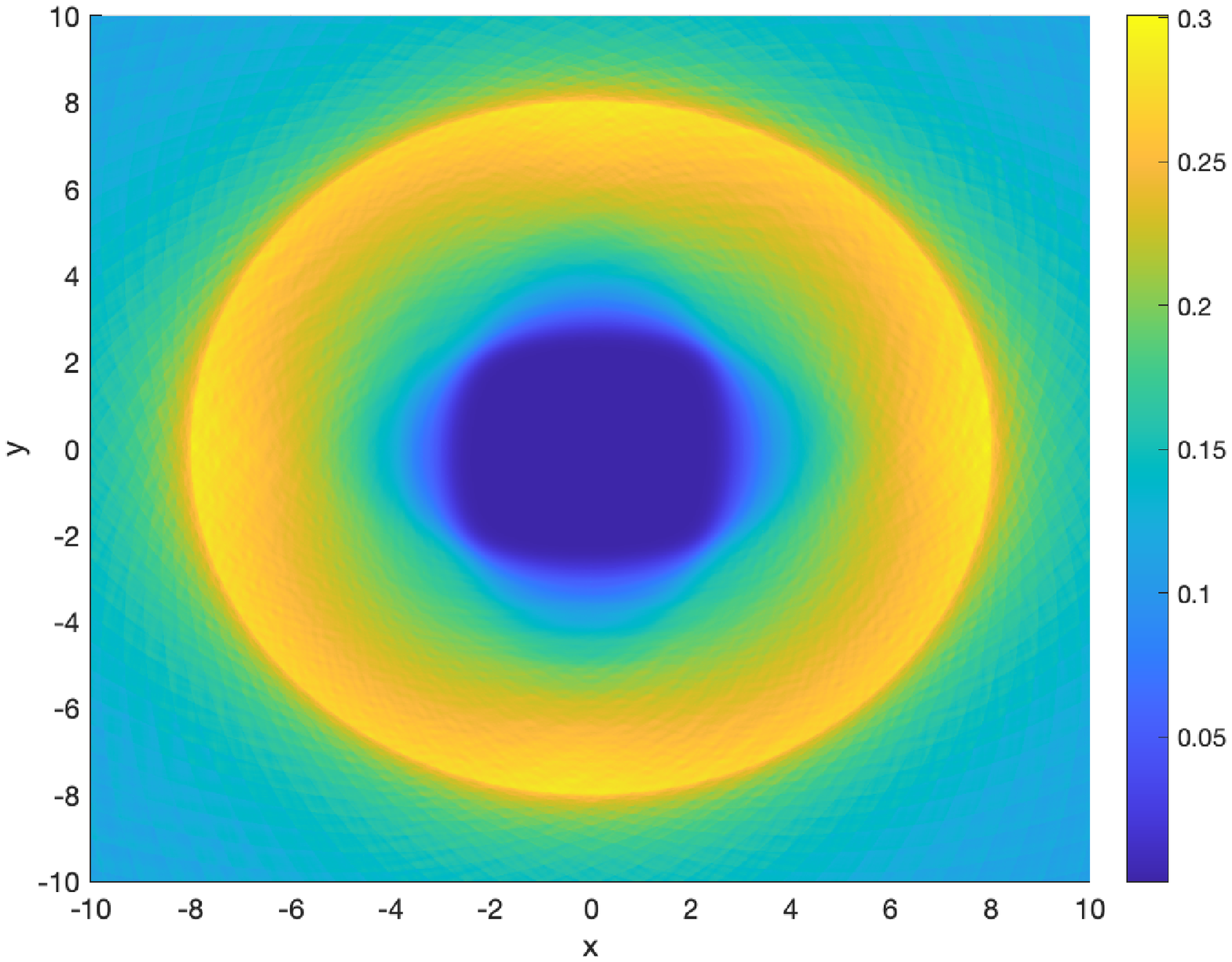}
		\caption*{$\theta=4$ rad}
	\end{minipage}
}
	\subfigure
{
	\begin{minipage}[t]{0.3\linewidth}
		\centering
		\includegraphics[scale=0.3]{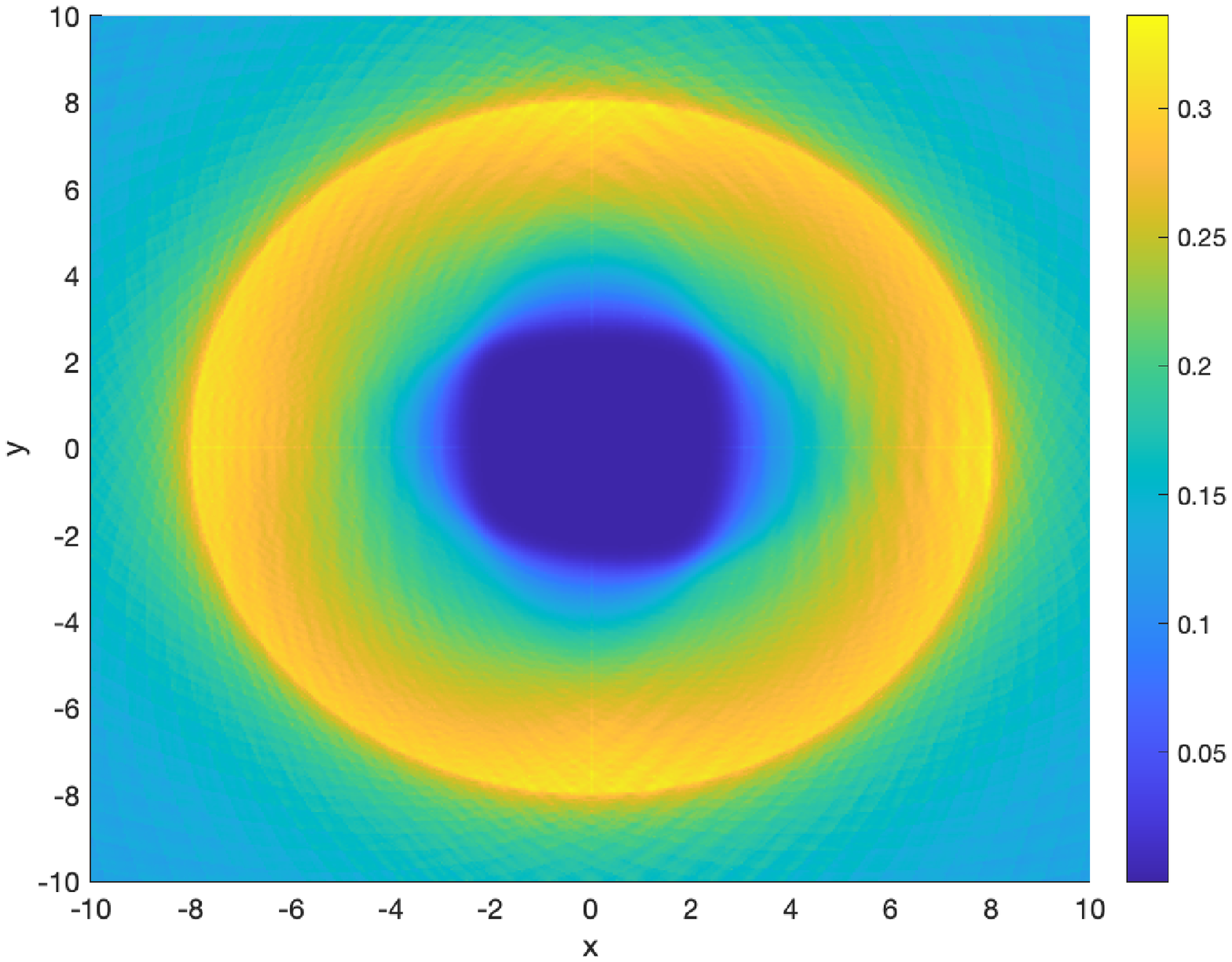}
		\caption*{$\theta=6$ rad}
	\end{minipage}
}
	\caption{Reconstruction of a sound-soft square $D=[-3,3]\times[-3,3]$ from one far-field pattern excited by a plane wave with different directions $d=(\cos\theta, \sin\theta)$.}
\label{Obstacle using different incident direction}
\end{figure}

Suppose that an incident plane wave $u^{i}(x)=e^{ikx\cdot d}$ with $k=6$ and $d=(\cos\theta,\sin\theta)$ is incident onto a sound-soft square $D=[-3,3]\times[-3,3]$. Set $R=8$ and $N=80$. To get an image of $D$, we use the spectral system of $F_{z,h}$ corresponding to sampling disks $B_h(z)$ of either the Dirichlet or impedance type. In the impedance case, the Robin coefficient on $\partial B_h(z)$ is set to be $\eta=-2+i$; see \eqref{imp-bound}. Recall that in the Dirichlet case, we have to assume that $k^2$ is not the Dirichlet eigenvalue of $-\Delta$ over $B_z(h)$ for all $h\in (0, 2R)$. The radii of these sampling disks are discretized as $h_m=2mR/M$ for $m=1,2,\cdots, M$. Numerics show that the choice of $M=160$ can avoid the Dirichlet eigenvalue problem.
In Fig. \ref{Obstacle using different boundary condition test disks}, we show the image of $D$ where the regularization parameter is $\alpha=1e-14$ and the incident angle is $\theta=4$ radians. We find that the image using test disks of impedance type is better than the image using sound-soft test disks. To compare the reconstruction results with different incident directions, we use test disks of impedance type and fix the regularization parameter at $\alpha=1e-14$. It is seen from
Fig. \ref{Obstacle using different incident direction} that the incident direction of a plane wave does not affect our imaging result too much.


In the aforementioned tests the noise level of the far-field pattern is set to be zero. Below we pollute the measurement data by
\ben
u_{D,\delta}^\infty(\hat x)=u_{D}^\infty(\hat x)+ \delta \kappa(\hat x) u_{D}^\infty(\hat x),
\enn
where $\kappa (\hat x)$ is a random function whose values are uniformly distributed between $-1$ to $1$ and $\delta>0$ is the noise level.
To test the influence of the noise, we fix $k=6$, $\theta=4$ rad, $R=8$, $M=160$, $N=80$, and $N_z=64$. We use test disks of sound-soft type. We again observe that a satisfactory image can be obtained only if the regularization parameter can be chosen properly. In our case, we take
$\alpha=1e-14$, $1e-11$ and $1e-4$ corresponding to the noise level $\delta=0$, $3\%$ and $8\%$.


\begin{figure}[H]
	\centering
	\subfigure
	{
		\begin{minipage}[t]{0.3\linewidth}
			\centering
			\includegraphics[scale=0.3]{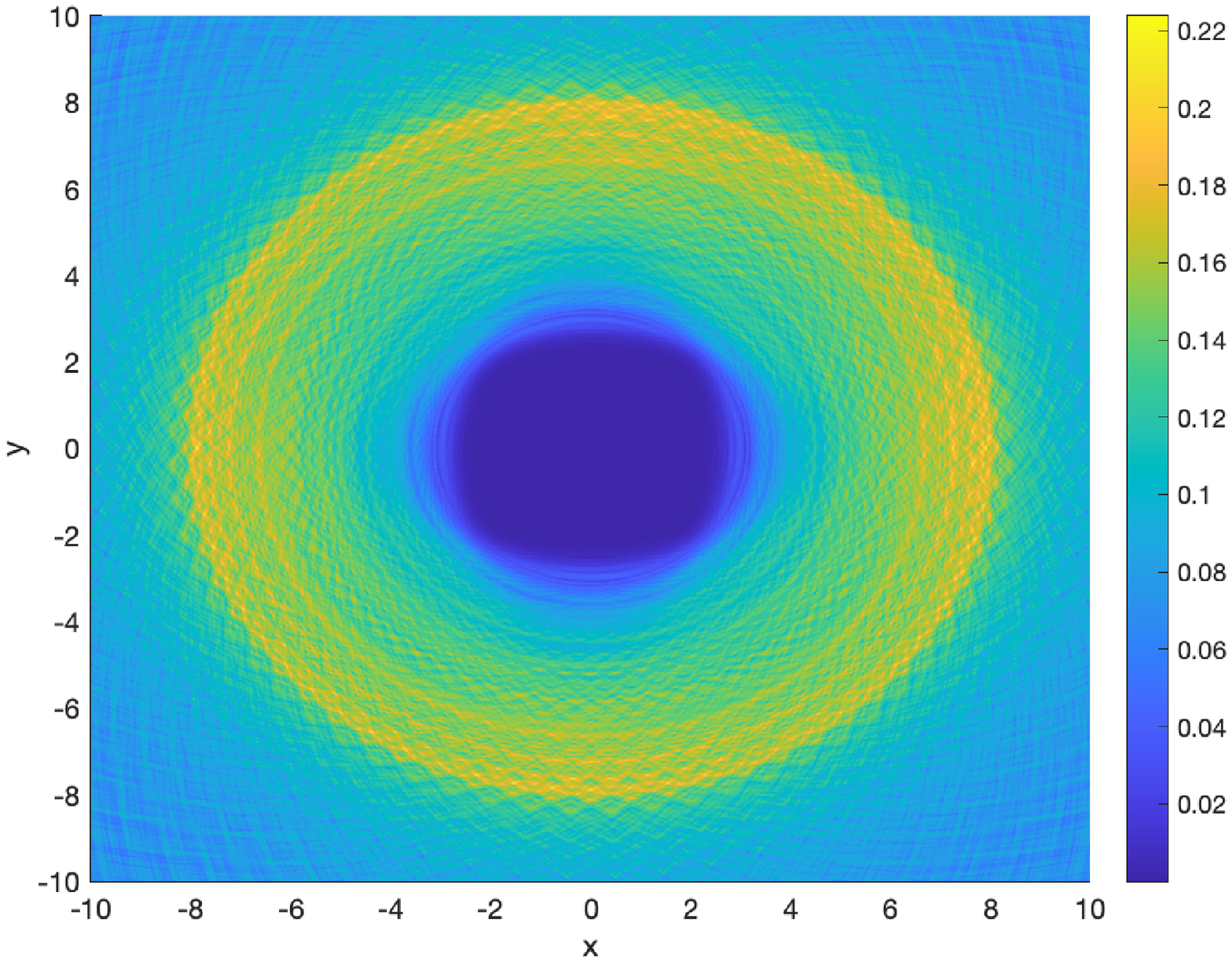}
			 \caption*{$\delta=0$}
		\end{minipage}
	}
	\subfigure
	{
		\begin{minipage}[t]{0.3\linewidth}
			\centering
			\includegraphics[scale=0.3]{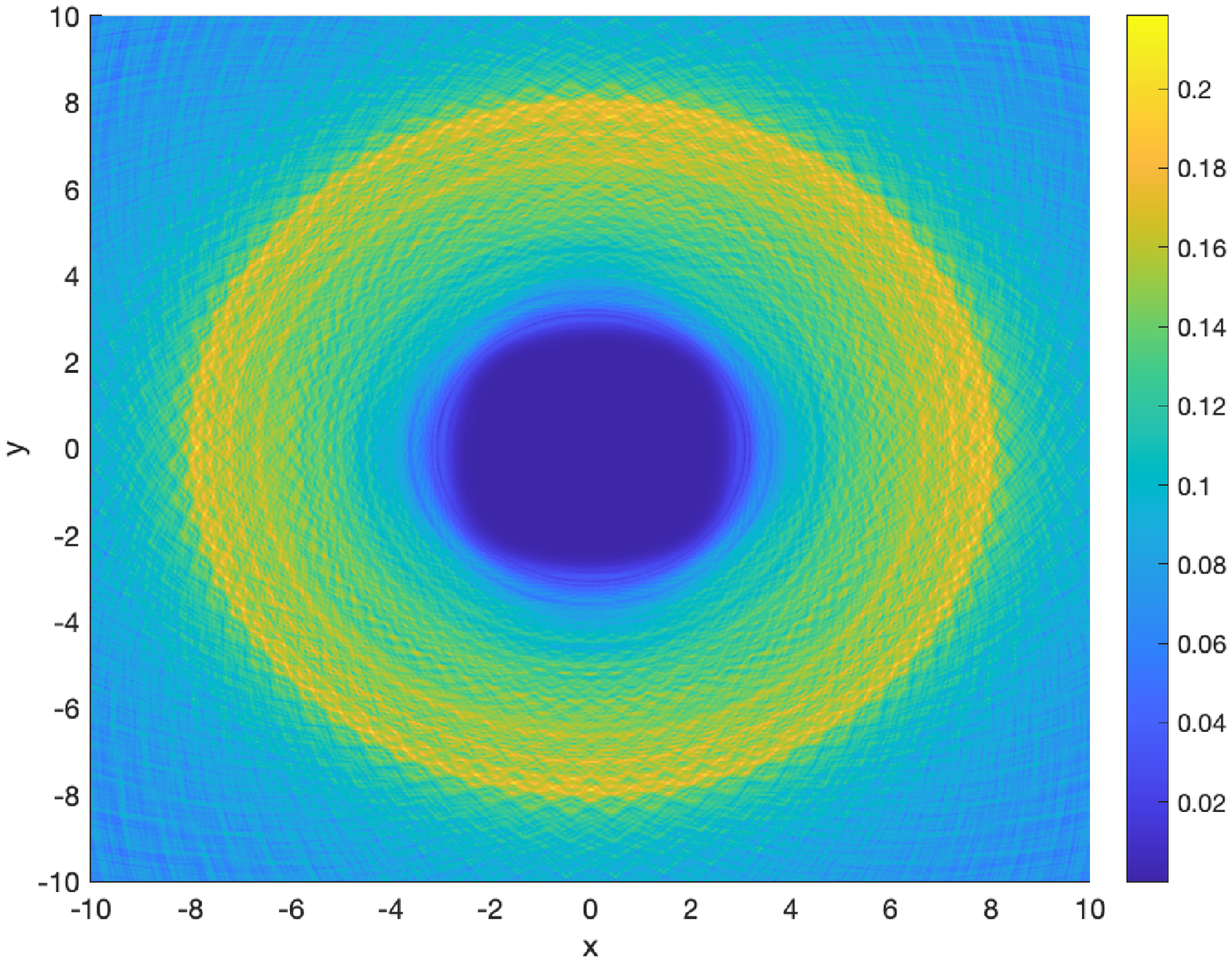}
			\caption*{$\delta=3\%$}
		\end{minipage}
	}
	\subfigure{
	\begin{minipage}[t]{0.3\linewidth}
		\centering
		\includegraphics[scale=0.3]{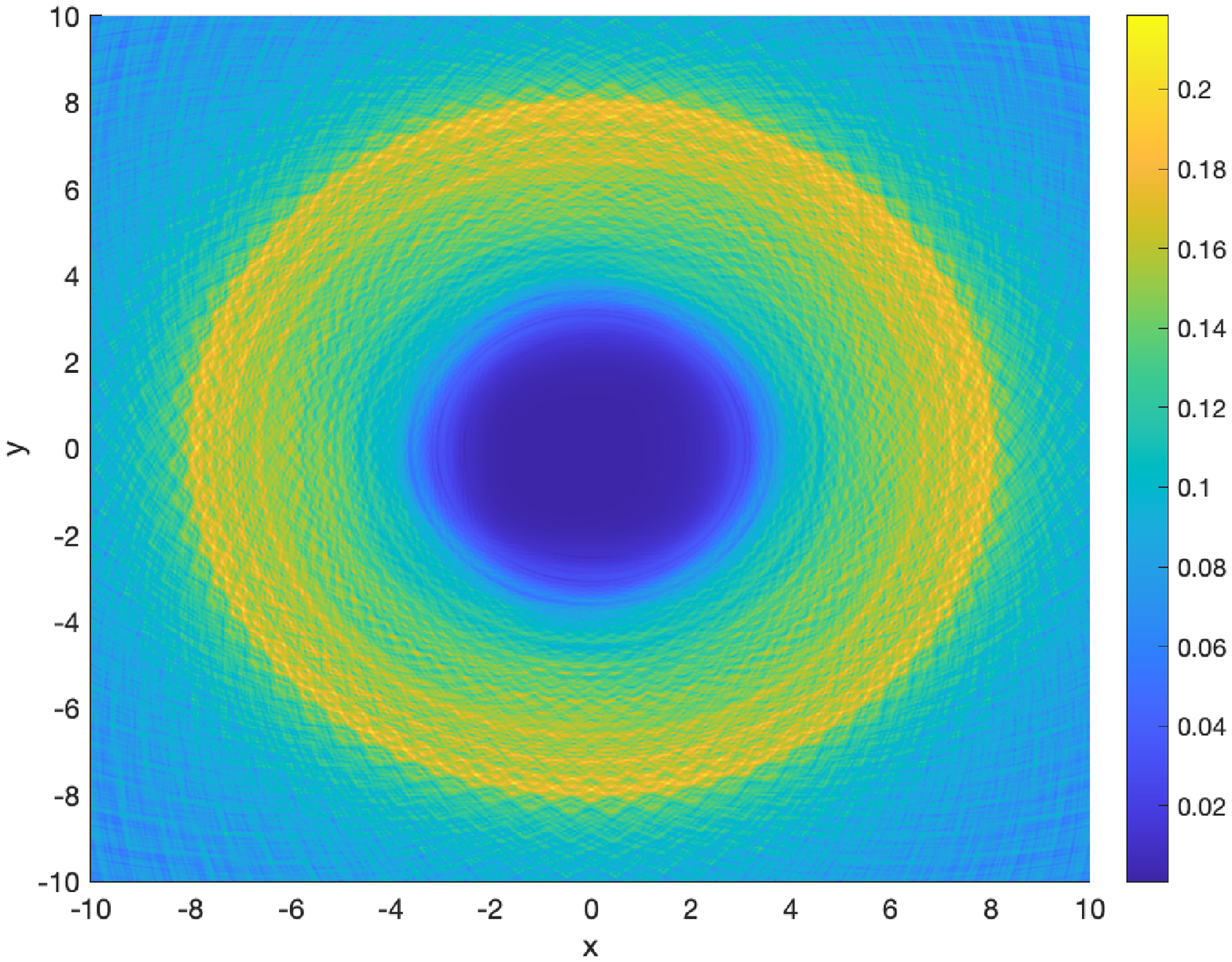}
		\caption*{$\delta=8\%$}
	\end{minipage}
}
	\caption{Reconstruction of a sound-soft square from polluted far-field pattern at the noise level $\delta$.}
    \label{Obstacle with different degrees of Gaussian noise}
\end{figure}

\section{Concluding discussions}\label{sec:5}
While the classical factorization method was motivated by the uniqueness proof in inverse scattering with infinitely many incident directions \cite{Isa1990,KK93}, its one-wave version, as explored here, originates from the unique determination of convex polygonal/polyhedral scatterers with a single incoming wave. The uniqueness proof implies that the wave field does not admit an analytical extension across a singular point lying on the interface. Being different from other domain-defined sampling methods, the imaging scheme proposed in this paper can be interpreted as a model-driven and data-driven approach, because it relies on both the Helmholtz equation and the a priori data for test scatterers.
In our numerics, these test scatterers are chosen as sound-soft or impedance disks, because the explicit forms of the spectra of the corresponding far-field operators have simplified numerical calculations. We remark that the a priori information of the unknown target $D$ can also be incorporated into the chosen test scatterers. Preliminary examples indeed show that the proposed scheme can be used to roughly capture the location and shape of a convex polygonal scatterer/source, due to the presence of interface singularities. The schemes developed in this paper can serve as the initial step for finding the boundary of an unknown target, when a single far-field pattern is available only. In the case of multi-static or multi-frequency measurement data, one can also design new imaging functionals based on (\ref{Wzh}).

Below we list several questions that are deserved to be further investigated in future.
\begin{itemize}
\item[(i)] Choice of the regularization parameter. Our numerical experience show that this parameter should at least depend on the noise level $\delta$, the radius $R$ of our test disks, the incident wave number $k$ as well as the truncation parameter $N$.
\item[(ii)] Convergence of the one-wave factorization method. Obviously, it is related to the blow-up rate of the function $h\rightarrow W(z,h)$ for $|z|=R$ as $|z-h|\rightarrow \max_{y\in D}|z-y|$. Subsection \ref{FMdisk} has shown that, for $D=\{z^*\}$,  we have $W(z,h)\sim -\ln (|z-z^*|-h)$ as $h\rightarrow |z-z^*|$. If $D$ is a polygonal scatterer, we conjecture that the convergence rate should rely on the singular behavior of the scattered field near corner points, and that a strongly/weakly singular corner could lead to a fast/slow convergence in detecting the corner.
 \item[(iii)] Further numerical tests by using other types of sample data $u^\infty(\hat x; \Omega)$. The test sample/scatterer $\Omega$ can be also sound-hard obstacles, penetrable scatterers and source terms, in addition to the sound-soft and impedance obstacles investigated in this paper. The far-field data $u^\infty(\hat x; \Omega)$ are also allowed to be excited at multi-frequencies. Hence, there is a variety of choices on the sample data $u^\infty(\hat x; \Omega)$ and on the shape and physical properties of $\Omega$.
 \item[(iv)] Analytical continuation test if $\partial D$ contains no singular points. Suppose that $D$ is a sound-soft obstacle with an analytic boundary. Our approach is capable of detecting some singular points of the analytical continuation of the wave field into the interior of $D$. It was shown in \cite{Millar} that these interior singularities of the Helmholtz equation are connected to the singularities of the Schwartz function of $\partial D$ (see \cite{Davis79}). Hence, information on $\partial D$ can still be extracted from a single far-field pattern, if we can disclose the relation between $\partial D$ and the singularities of the analytical continuation in $D$.
\end{itemize}

\section{Acknowledgements}
This work was supported by NSFC 11871092, NSFC 12071236 and NSAF U1930402.

\end{document}